\documentclass[12pt,a4paper]{amsart}

\usepackage[small, scriptlabels]{diagrams}
\newcommand{\bdi}{\begin{diagram}}
\newcommand{\edi}{\end{diagram}}

\usepackage{enumerate, amssymb}
\advance\hoffset-20mm \advance\textwidth40mm

\theoremstyle{plain}
\newtheorem{theorem}{Theorem}[section]
\newtheorem{thm}{Theorem}[section]
\newtheorem{cor}[thm]{Corollary}
\newtheorem{lem}[thm]{Lemma}
\newtheorem{prop}[thm]{Proposition}
\theoremstyle{definition}
\newtheorem{defi}[thm]{Definition}
\newtheorem{defis}[thm]{Definitions}
\newtheorem{conj}[thm]{Conjecture}
\newtheorem{conv}[thm]{Convention}
\newtheorem{nota}[thm]{Notation}
\newtheorem{rem}[thm]{Remark}
\newtheorem{rems}[thm]{Remarks}
\newtheorem{exa}[thm]{Example}
\newtheorem{exas}[thm]{Examples}
\newtheorem{prob}[thm]{Problem}
\newtheorem{probs}[thm]{Problems}
\newtheorem{ques}[thm]{Question}
\newtheorem{sit}[thm]{}

\newcommand{\Span}{ \operatorname{{\rm Span}}}
\newcommand{\Spec}{ \operatorname{{\rm Spec}}}

\newcommand{\sing}{ \operatorname{{\rm sing}}}

\newcommand{\Lie}{ \operatorname{{\rm Lie}}}

\newcommand{\Aut}{ \operatorname{{\rm Aut}}}
\newcommand{\Mat}{ \operatorname{{\rm Mat}}}
\newcommand{\LND}{ \operatorname{{\rm LND}}}
\newcommand{\GL}{ \operatorname{{\rm GL}}}
\newcommand{\SL}{ \operatorname{{\rm SL}}}
\newcommand{\ML}{ \operatorname{{\rm ML}}}
\newcommand{\FML}{ \operatorname{{\rm FML}}}
\newcommand{\Hom}{ \operatorname{{\rm Hom}}}
\newcommand{\im}{ \operatorname{{\rm im}}}
\newcommand{\trdeg}{ \operatorname{{\rm trdeg}}}

\newcommand{\C}{\ensuremath{\mathbb{C}}}
\newcommand{\T}{\ensuremath{\mathbb{T}}}

\newcommand{\Q}{\ensuremath{\mathbb{Q}}}
\newcommand{\Z}{\ensuremath{\mathbb{Z}}}
\newcommand{\N}{\ensuremath{\mathbb{N}}}
\newcommand{\G}{\ensuremath{\mathbb{G}}}

\newcommand{\kk}[1]{\bk^{[#1]}}

\newcommand{\fm}{{\mathfrak m}}

\newcommand{\sgoth}{{\ensuremath{\mathfrak{s}}}}
\newcommand{\lgoth}{{\ensuremath{\mathfrak{l}}}}

\newcommand{\cB}{{\ensuremath{\mathcal{B}}}}

\newcommand{\cG}{{\ensuremath{\mathcal{G}}}}

\newcommand{\cO}{{\ensuremath{\mathcal{O}}}}

\newcommand{\cH}{{\ensuremath{\mathcal{H}}}}
\newcommand{\cI}{{\ensuremath{\mathcal{I}}}}
\newcommand{\cM}{{\ensuremath{\mathcal{M}}}}

\newcommand{\cN}{{\ensuremath{\mathcal{N}}}}

\newcommand{\cY}{{\ensuremath{\mathcal{Y}}}}

\newcommand{\p}{\partial}
\newcommand{\id}{{\rm id}}
\newcommand{\vr}{\varrho}

\def\bals#1\eals{\begin{align*}#1\end{align*}}
\def\bal#1\eal{\begin{align}#1\end{align}}

\def\SAut{\mathop{\rm SAut}}
\def\kk{{\Bbbk}}
\def\AA{{\mathbb A}}

\def\kk{{\Bbbk}}

\def\PP{{\mathbb P}}
\def\V{{\mathbb V}}

\def\Der{\mathop{\rm Der}}

\def\deg{\mathop{\rm deg}}

\def\sl{\mathfrak{sl}}
\def\reg{{\mathop{\rm reg}}}

\def\rk{\mathop{\rm rk}}

\def\Sat{\mathop{\rm Sat}}

\def\codim{\mathop{\rm codim}}
\def\ML{\mathop{\rm ML}}
\def\Pic{\mathop{\rm Pic}}

\def\lto{\longrightarrow}
\def\ol{\overline}
\def\and{\quad\mbox{and}\quad}

\newcommand{\Cl}{ \operatorname{\rm Cl}}
\newcommand{\Pf}{ \operatorname{\rm Pf}}
\renewcommand{\div}{ \operatorname{\rm div}}

\renewcommand{\epsilon}{\varepsilon}
\renewcommand{\phi}{\varphi}

\newcommand{\bnum}{\begin{enumerate}}
\newcommand{\enum}{\end{enumerate}}
\renewcommand{\emptyset}{\varnothing}
\newcommand{\brem}{\begin{rem}}
\newcommand{\brems}{\begin{rems}}
\newcommand{\erem}{\end{rem}}
\newcommand{\erems}{\end{rems}}
\newcommand{\bprob}{\begin{prob}}
\newcommand{\eprob}{\end{prob}}
\newcommand{\bprobs}{\begin{probs}}
\newcommand{\eprobs}{\end{probs}}
\newcommand{\bques}{\begin{ques}}
\newcommand{\eques}{\end{ques}}
\newcommand{\bexa}{\begin{exa}}
\newcommand{\bexas}{\begin{exas}}
\newcommand{\eexa}{\end{exa}}
\newcommand{\eexas}{\end{exas}}
\newcommand{\bdefi}{\begin{defi}}
\newcommand{\edefi}{\end{defi}}
\newcommand{\bdefis}{\begin{defis}}
\newcommand{\edefis}{\end{defis}}
\newcommand{\bcor}{\begin{cor}}
\newcommand{\ecor}{\end{cor}}
\newcommand{\blem}{\begin{lem}}
\newcommand{\elem}{\end{lem}}
\newcommand{\bconv}{\begin{conv}}
\newcommand{\econv}{\end{conv}}
\newcommand{\bconj}{\begin{conj}}
\newcommand{\econj}{\end{conj}}
\newcommand{\bprop}{\begin{prop}}
\newcommand{\eprop}{\end{prop}}
\newcommand{\bthm}{\begin{thm}}
\newcommand{\ethm}{\end{thm}}
\newcommand{\bnota}{\begin{nota}}
\newcommand{\enota}{\end{nota}}
\newcommand{\bsit}{\begin{sit}}
\newcommand{\esit}{\end{sit}}
\newcommand{\be}{\begin{equation}}
\newcommand{\ee}{\end{equation}}
\newcommand{\bproof}{\begin{proof}}
\newcommand{\eproof}{\end{proof}}
\def\ba{\begin{array}}
\def\ea{\end{array}}

\thanks{
This work was done during a stay of the second, third, and fifth
authors at the Max Planck Institut f\"ur Mathematik at Bonn and a
stay of the first and the second authors at the Institut Fourier,
Grenoble. The authors thank these institutions for their
hospitality. The research of the forth author was partially
supported by Schweizerische Nationalfonds grant
$200020-134876/1$.}

\begin{document}
\title[Flexible varieties and
automorphism groups]
{Flexible varieties and automorphism groups}

\author{I.\ Arzhantsev, H.\ Flenner, S.\ Kaliman, F.\ Kutzschebauch,
M.\ Zaidenberg}
\address{Department of Algebra, Faculty of Mechanics and Mathematics,
Moscow State University, Leninskie Gory 1, GSP-1, Moscow, 119991,
Russia } \email{arjantse@mccme.ru}
\address{Fakult\"at f\"ur
Mathematik, Ruhr Universit\"at Bochum, Geb.\ NA 2/72,
Universit\"ats\-str.\ 150, 44780 Bochum, Germany}
\email{Hubert.Flenner@rub.de}
\address{Department of Mathematics,
University of Miami, Coral Gables, FL 33124, USA}
\email{kaliman@math.miami.edu}
\address{Mathematisches Institut,
Universit\"at Bern, Sidlerstrasse 5, CH-3012 Bern, Switzerland}
\email{frank.kutzschebauch@math.unibe.ch}
\address{Universit\'e Grenoble I, Institut Fourier, UMR 5582
CNRS-UJF, BP 74, 38402 St. Martin d'H\`eres c\'edex, France}
\email{Mikhail.Zaidenberg@ujf-grenoble.fr}

\begin{abstract} Given an irreducible affine algebraic variety $X$
of dimension $n\ge 2$, we let $\SAut (X)$ denote the special
automorphism group
 of $X$ i.e., the subgroup of the full automorphism group
$\Aut (X)$ generated by all one-parameter unipotent subgroups. We
show that if  $\SAut (X)$ is transitive on the smooth locus
$X_{\reg}$ then it is infinitely transitive on  $X_{\reg}$. In
turn, the transitivity is equivalent to the flexibility of $X$.
The latter means that for every smooth point $x\in X_{\reg}$ the
tangent space $T_xX$ is spanned by the velocity vectors at $x$ of
one-parameter unipotent subgroups of $\Aut (X)$. We provide also
various modifications and applications.

\end{abstract}

\maketitle
\date{today}

\thanks{
{\renewcommand{\thefootnote}{} \footnotetext{ 2010
\textit{Mathematics Subject Classification:}
14R20,\,32M17; secondary 14L30.\mbox{\hspace{11pt}}}}}

\vfuzz=2pt
\thanks{}

%{\footnotesize \tableofcontents}

\section*{Introduction}\label{intro}
All algebraic varieties and algebraic groups in this paper are
supposed to be reduced and defined over an algebraically
closed field $\kk$ of characteristic zero.
For such a variety $X$ we  say that a subgroup $H$ of the
automorphism group $\Aut(X)$ is {\em algebraic} if it admits a
structure of an algebraic group such that the natural map $H\times
X \to X$ is a morphism. We let $\SAut(X)$ denote the subgroup of
$\Aut (X)$ generated by all algebraic one-parameter unipotent
subgroups of $\Aut(X)$ i.e., algebraic subgroups isomorphic to the
additive group $\G_a$. The group $\SAut(X)$ is called the
{\em special automorphism group} of $X$; this is a normal subgroup
of $\Aut(X)$. In this paper we study transitivity properties of
the action of $\SAut(X)$ on an irreducible variety $X$.

For instance, the special automorphism group $\SAut (\AA^1)$
of the affine line is an algebraic group that acts transitively
but not 2-transitively. In contrast, for any $n\ge 2$ the group
$\SAut(\AA^n)$ is no longer  an algebraic group. Indeed, it
contains the infinite dimensional vector group of shears
$$
(x_1,\ldots,x_{n-1},y)\mapsto
(x_1,\ldots,x_{n-1},y+P(x_1,\ldots,x_{n-1}))\,,
$$
where $P\in \kk[x_1,\ldots,x_{n-1}]$ is an arbitrary
polynomial.
It is a well known and
elementary fact that $\SAut(\AA^n)$, $n\ge 2$, acts {\em
infinitely transitively} on $\AA^n$ that is, $m$-transitively for
any $m\ge 1$ (see e.g.\ \cite[Lemma 5.5]{KZ1} and references
therein).

There is a number of further cases, where $\SAut(X)$ acts
infinitely transitively. Consider, for instance, an equivariant
projective  embedding $Y\hookrightarrow\PP^n$ of a flag variety
$Y=G/P$. Then the special automorphism group of the affine cone
$X$ over $Y$ acts infinitely transitively on the smooth locus
$X_\reg$ of $X$ \cite[Theorem 1.1]{AKZ}. For non-degenerate toric
affine varieties of dimension $\ge 2$ a similar result is true
\cite[Theorem 2.1]{AKZ}. If $Y$ is an affine variety on which
$\SAut(Y)$ acts infinitely transitively, then the same holds for
the suspension \be\label{SuS} X=\{uv-f(y)=0\}\subseteq \AA^2\times
Y\ee over $Y$, where $f\in\cO(Y)$ is a non-constant function
(\cite[Theorem 3.2]{AKZ}; see also \cite[\S 5]{KZ1} for the case
$Y=\AA^n$).

Transitivity properties of the special automorphism group are
closely related to the flexibility of a variety, which was studied
in the algebraic context in \cite{AKZ}\footnote{In the analytic
context several other flexibility properties are surveyed in
\cite{For1}.}. The variety $X$ is called {\em flexible} if
every point $x\in X_\reg$ is. We say that a point $x\in X_\reg$ is
{\em flexible} if the tangent space $T_x X$ is spanned by the
tangent vectors to the orbits $H.x$ of one-parameter unipotent
subgroups $H\subseteq\Aut (X)$. Clearly, $X$ is flexible if one
point of $X_\reg$ is and the group $\Aut(X)$ acts transitively on
$X_{\rm reg}$.

The following theorem\footnote{Cf.\  Theorem \ref{maincor} below.}
confirms a conjecture formulated in \cite[\S 4]{AKZ}.

\begin{theorem}\label{mthm}
For an  irreducible affine variety $X$ of dimension $\ge 2$,
the following conditions are equivalent.
\begin{enumerate}[(i)]
\item The group $\SAut (X)$ acts transitively on  $X_\reg$.
\item The group $\SAut (X)$ acts
infinitely transitively on $X_\reg$.
\item  $X$ is a flexible variety.
\end{enumerate}
\end{theorem}

The varieties studied in \cite{AKZ} are flexible. As a
further example of a flexible variety one can consider the total
space of a homogeneous vector bundle over a flexible affine
variety (see Corollary \ref{5.23}). Every connected semisimple algebraic
group and, more generally, every connected linear algebraic group $G$
without non-trivial characters is generated by its unipotent
1-parameter subgroups, see \cite[Lemma 1.1]{Po2}. This implies
that any affine homogeneous space $G/H$ is flexible. Consequently
in case $\dim G/H\ge 2$ the special automorphism group
$\SAut(G/H)$ acts infinitely transitively. More generally, if a
semisimple algebraic group acts with an open orbit on a smooth
affine variety $X$ then $X$ is homogeneous with respect to a
bigger affine algebraic group without non-trivial characters and
so is flexible (see Theorem \ref{prohom}).

In contrast, a Lie group or an algebraic group $G$ cannot act
$m$-transitively on a variety $X$ for $m>\dim G/\dim X$.
Indeed, $m$-transitivity of the $G$-action on $X$ is equivalent
to the transitivity of the induced $G$-action on $X^m$ minus the
diagonals. According to A.\ Borel a much stronger fact is valid: a
real Lie group cannot act even $3$-transitively on a simply
connected, non-compact real manifold (see Theorems 5 and 6 in
\cite{Bo}). The latter remains true, without the assumption of
simple connectedness, for the actions of algebraic groups over
algebraically closed fields \cite[Korollar 2]{Kn};  cf.\
\cite{Po3} for related results.

Let us mention several applications. As an almost immediate
consequence it follows that  in a flexible irreducible
affine variety $X$ any finite subset $Z\subseteq X_\reg$ can be
interpolated  in $X_\reg$ by an $\AA^1$-curve, that is by a
curve isomorphic to $\AA^1$ and contained in $X_\reg$.
Indeed,  given a one-dimensional
$\G_a$-orbit $O$ in $X_{\reg}$ and a finite subset $Z'\subseteq O$
of the same cardinality as that of $Z$, by infinite transitivity
there is an automorphism $g\in\SAut (X)$ which sends $Z'$ to $Z$.
Then $g(O)\cong\AA^1$ is a $\G_a$-orbit passing through every
point of $Z$. In fact, $X_{\reg}$ is $\AA^1$-rich in the sense of
\cite{KZ2} (see Corollary \ref{5.41}). For the case $X=\AA^n$ this
is the Gromov-Winkelmann theorem, see \cite{Wi}.

An interesting class of  flexible varieties is formed by
degeneracy loci of generic matrices. These are the varieties
$X_r\subseteq\AA^{mn}$ consisting of $m\times n$-matrices of rank
$\le r$, where $1\le r\le\min(n,m)$. The elementary
transformations, which replace row $i$ by row $i+t\cdot$row $j$
($i\ne j, \, t\in \kk$), and similarly for columns, constitute
$\G_a$-actions on $X_r\setminus X_{r-1}$. By a standard fact of
linear algebra each matrix can be transformed to a normal form by
a sequence of elementary transformations.
Since $X_r\backslash X_{r-1}=(X_r)_\reg$  for $1\le r<\min (n,m)$,
in this range $\SAut(X_r)$ acts transitively and thus infinitely transitively
on $X_r\setminus X_{r-1}$ by Theorem \ref{mthm}.

We establish this infinite transitivity even
simultaneously for matrices of different ranks, see Theorem
\ref{transmat}. This shows that any finite collection of $m\times
n$ matrices can be diagonalized simultaneously by means of
elementary row- and column transformations depending polynomially
on the matrix entries. Similar statements also hold for symmetric
and skew-symmetric matrices, see Theorems \ref{transmatsym} and
\ref{transmatskew}. A related result on collective infinite
transitivity for conjugacy classes of matrices was established
earlier by Z.\ Reichstein \cite{Re1} using different methods.

The Gizatullin surfaces represent another interesting class of
examples, where flexibility manifests (see Example \ref{gisu}).
These are the normal affine surfaces which admit a
completion by a chain of smooth rational curves.

We provide as well a version of infinite transitivity involving
infinitesimally  near points. More precisely we prove\footnote{See
Theorem \ref{5.33}.}:

\setcounter{thm}{1}{}

\bthm\label{0.2} Let $X$ be a flexible irreducible affine
variety of dimension $n\ge 2$ equipped with an algebraic  volume
form\footnote{By this we mean a nowhere vanishing $n$-form defined
on $X_\reg$.} $\omega$. Then for every $m\ge 0$ and every finite
subset $Z\subseteq X_\reg$ there exists an automorphism $g \in
\SAut (X)$ with prescribed $m$-jets at the points $p\in Z$,
provided these jets preserve $\omega$ and inject $Z$ into
$X_\reg$. \ethm

In the analytic context similar results were obtained in
\cite{BF} and \cite{KK}.

Let us give a short overview of the content.

In Section \ref{fl-tr} we consider subgroups $G\subseteq\Aut(X)$
that are generated by collections of algebraic subgroups of
$\Aut(X)$ and discuss their general properties. Special classes of
such groups appeared  earlier on different occasions, see e.g.\
\cite[\S 5]{KZ1}, \cite{Po4}, and \cite{Po2}.
 Although $G$ is not in general an algebraic group,
 we show that the $G$-orbits
have the same properties as those of an algebraic group action
(cf.\ \cite{Ra}, Lemma 2). We give an extension of Kleiman's
Transversality Theorem in this context (see Theorem \ref{5.40})
and of the Rosenlicht Theorem on the separation of generic orbits
by rational invariants (see Theorem \ref{RT}). As an application
we  confirm a conjecture in \cite{Lie} concerning the field
Makar-Limanov invariant  (see \ref{fv.12.26.11}). We expect
further development of invariant theory for algebraically
generated groups.

 In Section \ref{inftrans} we  prove Theorem \ref{mthm}
(cf.\  Theorems \ref{maincor} and \ref{mainnew}). The
methods developed there are applied in Section
\ref{colinftr} to show infinite transitivity on several orbits
simultaneously, see Theorem \ref{maincorstrat}. This yields the
aforementioned application to matrix varieties.

Section \ref{solupro} contains the results on interpolation of
curves and automorphisms. In Section \ref{appli} we apply our
techniques to homogeneous spaces and their affine embeddings.

In the Appendix  we adopt a complex analytic point of view. We
show in particular that the Oka-Grauert-Gromov Principle is
available for smooth $G$-fibrations with flexible fibers, where
$G$ is an algebraically generated group of automorphisms (cf.\
Proposition \ref{6.10} and Corollary \ref{fds}).  Besides, we generalize the notion of flexibility to the analytic setting, and
compare our results with similar ones known in this
setup.

{\em Acknowledgment.} The referees made a number of pertinent
remarks that allowed us to  improve significantly the
presentation. We are grateful to all of them. Our thanks
are due also to Adrien Dubouloz and Marat Gizatullin for
useful observations.

\section{Flexibility versus transitivity}\label{fl-tr}
 In this section $X$ stands for a reduced and irreducible
algebraic variety.

\subsection{Algebraically generated groups of automorphisms}
Recall from the introduction that a subgroup $H$ of the
automorphism group $\Aut(X)$ is {\em algebraic} if $H$ has a
structure of an algebraic group such that the natural action $H
\times X \to X$ is a morphism. In the literature there are many attempts to define and to study more general classes of subgroups of $\Aut(X)$, see e.g.\ \cite{Ra}, \cite{Sha}\footnote{A thorough treatment of this approach can be found in \cite[Chapt.\ 4]{Ku}, along with some historical remarks and bibliography.}. For our purposes the following notion closely related to that in \cite[Definition 1.36]{Po2} turns out to be  useful.

\bdefi\label{agag}
A subgroup $G$ of $\Aut (X)$ is
called {\em algebraically generated} if it is generated as an
abstract group by a family $\cG$ of connected
algebraic\footnote{not necessarily affine.} subgroups of
$\Aut(X)$.
\edefi

The following notation will be useful in the sequel.

\bnota\label{1.0}
1. Let us introduce a partial order on the set of sequences in $\cG$
defined via
$$
(H_1,\ldots, H_m)\succcurlyeq (H'_1,\ldots H'_s)
\Longleftrightarrow
\exists i_1<\ldots <i_s:\quad (H_1',\ldots , H_s')=
(H_{i_1},\ldots, H_{i_s})\, .
$$
Clearly then any two sequences are dominated by a third one.

2. Given a sequence $\cH=(H_1,\ldots, H_s)$ in $\cG$ and a point $x\in X$
we consider
the  morphism \be\label{1.1.a} \Phi_{\cH,x}: H_1\times\ldots\times
H_s\to X, \qquad (h_1,\ldots, h_s)\mapsto (h_1\cdot\ldots\cdot
h_s).x\, .\ee \enota

\bprop\label{1.1} If the subgroup $G\subseteq\Aut(X)$ is
algebraically generated then for every point $x\in X$ the orbit
$G.x$ is locally closed. \eprop

\bproof Replacing $X$ by the Zariski closure of the orbit $G.x$ we
may assume that $X=\overline{G.x}$ i.e., the orbit of $x$ is dense
in $X$. Notice that for every finite sequence  $\cH=(H_1,\ldots,
H_s)$  in $\cG$ the subset $X_{\cH,x}=(H_1\cdot H_2. \cdot\ldots
\cdot H_s).x\subseteq X$ is constructible and irreducible, being
the image of the irreducible variety $H_1\times\ldots\times H_s$
under the morphism  $\Phi_{\cH,x}$. Taking a larger $\cH$ we
enlarge $X_{\cH,x}$ too (i.e.\ $X_{\cH,x} \subseteq X_{\cH',x}$
for $\cH' \succcurlyeq \cH$ ). By assumption the union of all such
sets $X_{\cH,x}$ is dense in $X$, hence also the union of the
closures $\bar X_{\cH,x}$ is. Since an increasing sequence of
closed irreducible subsets becomes stationary,  $X=\bar X_{\cH,x}$
for some $\cH$. In particular, the interior $\mathring X_{\cH,x}$
of $X_{\cH,x}$ is a non-empty open subset  in $X$. Now the
transitivity of the $G$-action on $G.x$ implies that $G.x$ is open
in $X$, as desired. \eproof

\brem\label{rem-L}  We are grateful to one  of the referees for pointing out to us that Proposition \ref{1.1} is contained in Lemma 2 of
\cite{Ra}.
\erem

For the next results it is useful to consider the map
\be\label{1.2.b} \Phi_\cH : H_1\times\ldots \times H_s\times X\lto
X\times X, \qquad (h_1,\ldots,h_s,x)\mapsto
(x,(h_1\cdot\ldots\cdot h_s).x)\,. \ee

\bprop\label{1.2} There are (not necessarily distinct) subgroups
$H_1,\ldots, H_s\in \cG$ such that \be\label{1.2.a} G.x= (H_1\cdot
H_2. \cdot\ldots \cdot H_s).x\quad \forall x\in X. \ee \eprop

\bproof The image $Z_\cH=\Phi_\cH(H_1\times\ldots \times H_s\times
X)$ in $X\times X$ is constructible and irreducible. As before, if
$\cH' \succcurlyeq \cH$ then $Z_{\cH} \subseteq Z_{\cH'}$. Hence
the union  of closures $Z=\bigcup_\cH{ \bar Z_\cH}$ stabilizes in
$X\times X$ i.e., it coincides with $\bar Z_\cH$ for $\cH$
sufficiently large. In particular $Z$ is closed.

Let $\mathring Z_\cH$ be the interior of $Z_\cH$ in $Z$. It
follows as before that also $\{\mathring Z_\cH\}$ becomes
stationary and that the union  $Z'=\bigcup_\cH \mathring Z_\cH$ is
an open dense subset of $Z$.

Consider the $G$-action on $X\times X$ given by
$g.(x,y)=(g.x,y)$. If $\cH=(H_1,\ldots, H_s)$ and $H\in \cG$ then
for any $(h_1,\ldots, h_s)\in H_1\times\ldots\times H_s$ and $h\in
H$ we have
$$
h.\Phi_{\cH}(h_1,\ldots, h_s,x)=
h.(x,(h_1\cdot\ldots\cdot h_s).x)
=\Phi_{(\cH,H)}(h_1,\ldots, h_s,h^{-1},h.x).
$$
Hence $h.Z_\cH\subseteq Z_{(\cH,H)}$.
It follows that $Z$ and $Z'$ are $G$-invariant.

Consider now for $\cH$ sufficiently large the sets $Z_\cH$,
$Z'=\mathring Z_\cH$, and $Z=\bar Z_\cH$ as families
over $X$ via
the first projection $p: (x,y)\mapsto x$. By \cite[9.5.3]{EGA}
there is an open dense subset $V$ of $X$ such that $Z'(x)$ is
dense in $Z(x)$ for all $x\in V$, where for a subset $M\subseteq
X\times X$ we denote by $M(x)$ the fiber of $p|M:M\to X$ over $x$.
Since $Z$ and $Z'$ are invariant under the action of $G$ and the
projection $p$ is equivariant, we may suppose that $V$ is as well
$G$-invariant.

In particular there is a sequence $\cH_0$ such that
$Z_\cH(x)=(H_1\cdot\ldots\cdot H_s).x$ is dense in $Z(x)$ for all
$x\in V$ and all sequences $\cH=(H_1,\ldots, H_s)$ dominating
$\cH_0$. It follows that $Z(x)$ is closure of  the orbit $G.x$
and so $(H_1\cdot\ldots\cdot H_s).x$ is dense in the orbit
$G.x$ for all $x\in V$.

We claim that for every point $x\in V$
$$
(H_s\cdot\ldots\cdot H_1\cdot H_1\cdot\ldots\cdot  H_s).x=G.x\,.
$$
Indeed, for any $y\in G.x$ the sets $(H_1\cdots H_s).x$ and
$(H_1\cdots H_s).y$ are both dense in the orbit $G.x=G.y$. Hence
they have a common point, say $z$. Thus $$y\in (H_s
\cdot\ldots\cdot H_1).z\subseteq (H_s\cdot\ldots\cdot H_1\cdot
H_1\cdot\ldots\cdot H_s).x\,.$$ Replacing $\cH$ by the larger
sequence $(H_s,\ldots ,H_1,H_1,\ldots, H_s)$ it follows  that
$$
(H_1\cdot\ldots\cdot H_s).x=G.x \mbox{ for all $x\in V$
simultaneously. }
$$
The complement $Y=X\backslash V$ is closed, $G$-invariant,
and all its irreducible components are of dimension $<\dim X$.
Using induction on the dimension of $X$ we see
that \eqref{1.2.a}
holds for $\cH$ sufficiently large
and all $x\in X$ simultaneously,
concluding the proof.
\eproof

\brem \label{char} Propositions \ref{1.1} and
\ref{1.2} remain true with the same proofs for varieties over
algebraically closed fields of arbitrary characteristic.

2. \label{orbitsub} In the setup of Proposition \ref{1.2}, if
$Y\subseteq X$ is constructible then so is its `orbit' $G.Y$.
Indeed, by Proposition \ref{1.2} for some $\cH = (H_1, \ldots ,
H_s)$ this orbit is the image of $H_1 \times \ldots \times H_s
\times Y$ under the composition of $\Phi_\cH$ and the natural
projection $X \times X \to X$ to the second factor.
\erem

\bdefi\label{1.3} A sequence $\cH=(H_1,\ldots , H_s)$ in $\cG$
satisfying condition \eqref{1.2.a} of \ref{1.2} will be called
{\em complete}. \edefi

%\brem \erem
%We note that this not any longer true for the property  `locally
%closed' even in the case of algebraic group actions (e.g.,
%consider the orbit of the $x$-axis in $\C^2_{x,y}$ under the
%$\G_a$-action $t.(x,y) = (x, y +tx)$). \erem

\bprop\label{1.5} Assume that the generating family $\cG$ of
connected algebraic subgroups is closed under conjugation in $G$,
i.e., $gHg^{-1}\in \cG$ for all $g\in G$ and $H\in \cG$. Then
there is a sequence $\cH=(H_1,\ldots, H_s) $ in  $\cG$ such that
for all $x\in X$ the tangent space $T_x(G.x)$ of the orbit $G.x$
is spanned by the tangent spaces
$$
T_x(H_1.x),\ldots ,T_x(H_s.x)\,.
$$
\eprop

\bproof We claim that $T_x(G.x)$ is spanned by the tangent spaces
$T_x(H.x)$, where $H\in \cG$. Indeed, consider a complete sequence
$H_1,\ldots, H_s\in \cG$ such that the map $\Phi_{\cH,x}:H_1\times
\ldots \times H_s\to G.x$ in \eqref{1.1.a} is surjective. By
\cite{Ha}[III, Corollary 10.7] this map is generically smooth.
Thus for some point $y=(h_1\cdot\ldots\cdot h_s).x\in G.x$ the
tangent map
$$
d\Phi_{\cH,x}: T_{(h_1, \ldots, h_s)}
(H_1\times \ldots \times H_s)\lto T_y(G.x)
$$
is surjective. Multiplication by $g=(h_1\cdot\ldots\cdot
h_s)^{-1}$ yields an isomorphism $\mu_g:G.x\to G.x$ which sends
$y$ to $x$. Hence the composition $\mu_g\circ\Phi_{\cH,x}$ has a
surjective tangent map

\be\label{eq1.8}
d(\mu_g\circ \Phi_{\cH,x}): T_{(h_1, \ldots, h_s)}
(H_1\times \ldots \times H_s)\cong
\prod_{\sigma=1}^s T_{h_\sigma}(H_\sigma)\lto T_x(G.x)\,.
\ee
Letting $g_\sigma:= h_{\sigma +1}\cdot\ldots\cdot h_s$ the isomorphism
$$
\tilde H_\sigma:= h_1\times \ldots \times h_{\sigma-1}\times H_\sigma
\times h_{\sigma +1}\times\ldots \times h_s\lto
H_\sigma':=g_\sigma ^{-1}H_\sigma g_\sigma
$$
with $( h_1, \ldots , h_{\sigma-1}, h, h_{\sigma +1},\ldots , h_s)
\mapsto g_\sigma^{-1}h_\sigma^{-1} hg_\sigma$
identifies the restriction $\mu_g\circ \Phi_{\cH,x}|\tilde H_\sigma$ with
$$
\varphi: H_\sigma'=g_\sigma ^{-1}H_\sigma g_\sigma\to G.x,\qquad
h'\mapsto h'.x\,.
$$
Thus the restriction of the tangent map $d(\mu_g\circ
\Phi_{\cH,x})$ to the factor $T_{h_\sigma}(H_\sigma)$ in
\eqref{eq1.8} can be identified with the map $d_e\varphi:
T_e(H'_\sigma)\to T_x(G.x)$, where $e$ denotes the identity
element of $H'_\sigma$. Thus the claim follows.

Consider further the map $\Phi_\cH:H_1\times\ldots \times
H_s\times X\to Z\subseteq X\times X$ as in \eqref{1.2.b}
associated with a complete sequence $\cH$, where $Z\subseteq
X\times X$ is the closure of the image of $\Phi_\cH$. Choose an
invariant open subset $V\subseteq X_\reg$ such that the first
projection $p: Z\to X$ is smooth over $V$. Note that the fiber of
$Z_V=p^{-1}(V)\to V$ over $x$ is just the orbit $G.x$. Let us
consider the map of relative tangent bundles
$$
d\Phi_\cH: T(H_1\times \ldots \times H_s\times V/V)\to
\Phi_\cH^*(T(Z_V/V))\,
$$
and its restriction to $(e,\ldots,e)\times V\cong V$, where $e$
is the identity element in $G$ and therefore in each $H_i$,
$$
d\Phi_\cH: T_eH_1\times \ldots \times T_eH_s
\times V\to \Phi_\cH^*(T(Z_V/V))|V\,.
$$
The set $U_\cH$ of points in $V$ where this map is surjective, is
open. By the above claim, the union $\bigcup_\cH U_\cH$ coincides
with $V$. Any two sequences $\cH_1$ and $\cH_2$ are dominated by a
third $\cH_3$ in the partial order as in \ref{1.0}, and the
corresponding subset $U_{\cH_3}$ contains both $U_{\cH_1}$ and
$U_{\cH_2}$. Thus the increasing union $\bigcup_\cH U_\cH$
stabilizes, that is, $V=U_\cH$ for $\cH$ sufficiently large.
Induction on the dimension of $X$ as in the proof of Proposition
\ref{1.2} ends the proof. \eproof

\brem\label{1.65} It may happen for a family $\cG$
which is not closed under conjugation
that  for some point $x\in X$ the tangent spaces
$$
T_x(H_1.x), \ldots, T_x(H_s.x)
$$
do not span $T_x(G.x)$, whatever is the sequence $\cH=(H_1,\ldots,
H_s)$ in $\cG$. For instance, the group $G=\SL_2$ is generated by
the family $\cG=\{U^+,U^-\}$, where $U^\pm$ are the subgroups of
upper and lower triangular unipotent matrices. Letting $\SL_2$ act
on itself by left multiplication the tangent space $T_eG$ of
the orbit $G=G.e$ is $\sl_2$, while for any sequence
$\cH=(H_1,\ldots,H_s)$ in $\cG$ the tangent spaces $T_e(H_1),
\ldots ,T_e( H_s)$ are  contained in the $2$-dimensional subspace
$T_e(U^+)+T_e(U^-)$.
\erem

\bdefi\label{gflex} Let $G\subseteq \Aut(X)$ be algebraically
generated by a family $\cG$ of connected algebraic subgroups,
which is closed under conjugation. We say that a point $p\in
X_{\reg}$ is {\em $G$-flexible} if the tangent space $T_p X$ at
$p$ is generated by the subspaces $T_p(H.p)$, where $H\in \cG$.
\edefi

\bcor\label{gflex1} With $G$ and $\cG$ as in Definition \ref{gflex}
the following hold.

\bnum[(a)]
\item A point $p\in X_\reg $ is $G$-flexible if and
only if the orbit $G.p$ is open in X.
\item An open $G$-orbit (if it
exists) is unique and consists of all $G$-flexible points in $X_\reg$.
\enum
\ecor

\bproof  (a) By Proposition \ref{1.5} the morphism $\Phi_{\cH ,
p}:H_1\times\ldots\times H_s\to G.p$ is surjective and smooth for
an appropriate choice of a sequence $\cH$. Now (a) follows.
Furthermore (b) follows from (a) since any two open $G$-orbits
overlap and so must coincide. \eproof

Let us note that by Corollary \ref{gflex1}(a) the
definition of a $G$-flexible point only depends on $G$ and not on
the choice of the generating set $\cG$.

Using the semicontinuity of the fiber dimension we can deduce the
following semicontinuity result for the orbits of algebraically
generated groups.

\bcor\label{1.4} If a group $G\subseteq \Aut (X)$ is algebraically
generated then the function $x\mapsto \dim G.x$ is lower
semicontinuous on $X$. In particular, there is a Zariski open
subset $U\subseteq X$ filled in by orbits of maximal dimension.
\ecor

\bproof We may suppose that $G=\langle\cG\rangle$, where $\cG$ is
a family of connected algebraic subgroups of $\Aut (X)$ closed
under conjugation in $G$. For a  complete sequence
$\cH=(H_1, \ldots, H_s)$  consider the map $\Phi_\cH$ from
\eqref{1.2.b}. By the semicontinuity of fiber dimension the
function
$$
X\ni x\longmapsto \dim_{\tau(x)} \Phi_\cH^{-1}(x,x)
$$
is upper semicontinuous on $X$, where $\tau(x)=(1,\ldots, 1,x)
\in H_1 \times \ldots \times H_s\times X$. Here $\Phi_\cH^{-1}(x,x)$
is just the fiber of the map
$\Phi_{\cH,x}: H_1\times\ldots\times H_s\to G.x$ over $x$.

Fix a point $x_0\in X$. Enlarging $\cH$ we may assume that
$\Phi_{\cH,x_0}$ is a submersion. Thus for $x$ in a suitable
neighborhood $U$ of $x_0$ \bals
\dim G.x_0=&\sum_{\sigma=1}^s\dim H_\sigma -\dim \Phi^{-1}(x_0,x_0)\\
\le &\sum_{\sigma=1}^s\dim H_\sigma -\dim \Phi^{-1}(x,x)\\
\le & \dim G.x. \eals It follows that $\dim G.x\ge \dim G.x_0$ for
$x\in U$, as required. \eproof

The following analog of the Rosenlicht Theorem on rational
invariants holds for algebraically generated groups. The proof of
this theorem given in \cite[Theorem 2.3]{PV} works {\em mutatis
mutandis} in our setting due to Proposition \ref{1.2} and
Corollary \ref{1.4}.

%
%In view of our preceding results the following analog of the
%Rosenlicht Theorem on rational invariants holds in our setting
%with an almost identical proof, see e.g.\ \cite[Theorem 2.3]{PV}.
%For the reader's convenience we add the argument.

\bthm\label{RT} Let $G$ be an algebraically generated group acting
on $X$. Then there exists a finite collection of
rational $G$-invariants which separate $G$-orbits in general
position. \ethm

\bproof 
Replacing $X$ by a subset $U$ as in Corollary \ref{1.4} we
may assume that all orbits of $G$ are of maximal dimension. In
particular then all $G$-orbits are closed in $X$.  Let
$\Gamma\subseteq X \times X$ consist of all pairs $(x,x')$ such
that $x$ and $x'$ are in the same $G$-orbit. Note that this is
just the image of the map $G\times X\to X\times X$ with
$(g,x)\mapsto (g.x,x)$. As we have seen in the proof of
Proposition \ref{1.2}, $\Gamma$ contains an open dense subset, say
$\Gamma_0$, of the closure $\overline\Gamma$ in $X \times X$.

Letting $G$ act on the first component of $X\times X$ we may
assume that $\Gamma_0$ is $G$-invariant, since otherwise we can
replace it by the union of all translates of $\Gamma_0$. If
$p_2:\Gamma_0\to X$ denotes the second projection then for a
general point $x\in X$ the fibre $p_2^{-1}(x)=G.x\times\{x\}$ is
closed in $X\times \{x\}$. Hence there is an open dense subset
$U\subseteq X$ such that $\Gamma_0\cap p_2^{-1}(U)$ is closed in
$X\times U$. In particular it follows that $\Gamma\cap X\times U$
is closed in $X\times U$. Shrinking $U$ we may also assume that
$U$ is affine.

Let $\cI\subseteq \cO({X\times U})$ be the ideal sheaf of
$\Gamma\cap X\times U$, and let $J$ be the ideal generated
by $\cI$ in the algebra $\cM er(X) \otimes \cO
(U)$\footnote{Here $\cM er(X)$ denotes the function field of
$X$.}. The ideal $J$ is $G$-invariant assuming that $G$ acts on
the first factor of $X\times X$. Moreover, $J$ is generated as a
 $\cM er(X)$-vector subspace by $G$-invariant elements (see
\cite[Lemma 2.4]{PV}). We can find a finite set of generators of
$J$, say $F_1,...,F_p$, among these elements. We have
$$F_i = \sum_s f_{is}\otimes u_{is},\quad\mbox{where}\quad f_{is}
\in \cM er(X)^G\quad\mbox{and}\quad u_{is} \in \cO(U)\,.$$
Let us show that the functions $f_{is}$ separate orbits in general
position.

Shrinking $U$ once again we may assume that all the $f_{is}$ are
regular functions on $U$ and that the elements $F_i$ generate the
ideal $\cI$. Then the orbit of a point $x\in U$ is defined by the
equations $F_i(x,y)=\sum_s f_{is}(x) u_{is}(y)=0$, $i=1,...,p$.
Consequently, the equalities $f_{is}(x_1)=f_{is}(x_2)$ for all $i$
and $s$ imply that $G.x_1=G.x_2$ on $U$. \eproof

As in \cite[Corollary on p.\ 156]{PV} this theorem
has the following consequence.

\bcor\label{RTcor} Let $G$ be an algebraically generated group
acting on $X$.
% a variety $X$ of dimension $n$.
Then
$$
\trdeg (\cM er(X)^G:\kk)=\min_{x\in X}\,\{ {\codim}_{X} G.x\}\,,
$$
where $\cM er(X)$ stands for the rational function field of
$X$. In particular, $G$ has an open orbit in $X$  if and only
if $\cM er(X)^G=\kk$. \ecor

\subsection{Transversality}\label{trthm}

If a connected algebraic group $G$ acts transitively on an
algebraic variety $X$ and  $Y$, $Z$ are smooth subvarieties of $X$
then by Kleiman's Transversality Theorem \cite{Klei} a general
$g$-translate $g.Z$ ($g\in G$) meets $Y$ transversally. In this
subsection we extend Kleiman's Theorem to the case of an arbitrary
algebraically generated group.

\bthm\label{5.40} Let a subgroup $G\subseteq \Aut(X)$ be
algebraically generated by a system $\cG$ of connected algebraic
subgroups closed under conjugation in $G$. Suppose that $G$ acts
with an orbit $O$ open in $X$.

Then there exist subgroups $H_1,\ldots, H_s\in \cG$ with the
following property: \\
For any locally closed reduced subschemes $Y$ and $Z$ in $O$
one can find a Zariski dense open subset $U=U(Y,Z)\subseteq
H_1\times \ldots \times H_s$ such that every element $(h_1,\ldots,
h_s)\in U$ satisfies the following conditions.
\bnum[(a)]
\item  The translate $(h_1\cdot\ldots\cdot h_s).Z_\reg$
meets $Y_\reg$ transversally.

\item $\dim (Y\cap (h_1\cdot\ldots\cdot h_s).Z)\le
\dim Y+\dim Z-\dim X$ \footnote{We let the dimension of the
empty set be equal to $-\infty$.}.\\ In particular $Y\cap
(h_1\cdot\ldots\cdot h_s).Z=\emptyset$ if $\dim Y+\dim Z<\dim X$.
\enum \ethm

The proof is based on the following auxiliary result, which is
complementary to Proposition \ref{1.5}.

\bprop\label{auxi} Let the assumption of Theorem \ref{5.40} hold.
Then there is a sequence $\cH=(H_1,\ldots, H_s)$ in $\cG$ so that
for a suitable open dense subset $U\subseteq H_{s}\times\ldots
\times H_{1}$ \footnote{The inverse enumeration here is
convenient for applying recursion.} the map \be\label{mapop}
\Phi_{s}: H_s\times \ldots\times H_1\times O\lto O\times O
\quad\mbox{with} \quad (h_s,\ldots,h_1,x)\mapsto
((h_s\cdot\ldots\cdot h_1).x ,x) \ee is smooth on $U\times O$.
\eprop

\bproof According to Proposition \ref{1.2}  there are subgroups
$H_1,\ldots, H_s\subseteq G$  in $\cG$ such that
$\Phi_{s}$ is surjective. Hence there is an open dense subset
$U_s\subseteq
H_s\times \ldots\times H_1\times O$ on which $\Phi_s$ is smooth.
Assuming that $U_s$ is maximal with this property we consider the
complement $A_s=(H_s\times \ldots\times H_1\times O)\backslash
U_s$.

Let us study
the effect of increasing the number of factors, i.e., passing to
$$
\Phi_{s+1}: H_{s+1}\times \ldots\times H_1\times O\lto O\times O
$$
The map $\Phi_{s+1}$ is smooth on $H_{s+1}\times U_s$. Indeed,
for every $h_{s+1}\in H_{s+1}$
we have a commutative diagram
\bdi[small]
H_{s+1}\times \ldots\times H_1\times O &
\rTo^{\Phi_{s+1}} & O\times O\\
\uTo<{h_{s+1}\times \id }&&   \uTo>{h_{s+1}\times\id}\\
\{1\}\times H_s\times \ldots\times H_1\times O&\rTo^{\Phi_s} &
O\times O \edi where the lower horizontal map is smooth on  $U_s$.
In other words, $U_{s+1}\supseteq H_{s+1}\times U_s$ or,
equivalently, $A_{s+1}\subseteq H_{s+1}\times A_s$. We claim that
increasing the number of factors by $H_{s+1},\ldots, H_{s+t}$ in a
suitable way, we can achieve that
\be\label{auxia} \dim A_{s+t}<
\dim (H_{s+t}\times\ldots\times H_{s+1}\times A_s) \,.
\ee
If $(h_s,\ldots,h_1,x)\in A_s$ and $y=(h_s\cdot\ldots\cdot h_1).x$
then for suitable $H_{s+t},\ldots, H_{s+1}$ the map
$$
H_{s+t}\times\ldots\times H_{s+1}\times O\lto O\times O
\quad\mbox{with} \quad (h_{s+t},\ldots,h_{s+1},x)\mapsto
((h_{s+t}\cdot\ldots\cdot h_{s+1}).x ,x)
$$
is smooth at all points $(e,\ldots, e, y)$, where $e$ is the
identity element of $G$; see Proposition \ref{1.5}. In particular
$\Phi_{s+t}$ is smooth at all points $(e,\ldots, e, h_s,\ldots,
h_1,x)$ with $x\in O$, i.e.\
$$
(e,\ldots, e)\times A_s\cap  A_{s+t}=\emptyset.
$$
Now \eqref{auxia} follows.

Thus increasing the number of factors suitably we can achieve
that\footnote{In fact we can make the difference
$\dim(H_s\times\ldots\times H_1)-\dim A_s$ arbitrarily large.}
$$\dim A_s<\dim(H_s\times\ldots\times H_1)\,.$$
In particular, the image of $A_s$ under the projection
$$
\pi: H_{s}\times\ldots\times H_{1}\times O\lto
H_{s}\times\ldots \times H_{1}
$$
is contained in a proper, closed subvariety of
$H_{s}\times\ldots \times H_{1}$. Hence there is an open dense
subset $U\subseteq H_{s}\times\ldots \times H_{1}$ such that
$\Phi_s: U\times O\to O\times O$ is smooth. \eproof

\bproof[Proof of Theorem \ref{5.40}] Let us first show (a).
Replacing $Y$ and $Z$ by $Y_\reg$ and $Z_\reg$, respectively, we
may assume that $Y$ and $Z$ are smooth.  Applying  Proposition
\ref{auxi}  there are subgroups $H_1,\ldots, H_s$ in $\cG$ such
that $\Phi_s:U \times O\to O\times O$ is smooth for some open
subset $U\subseteq H_1\times\ldots \times H_s$. In particular
$\cY=\Phi_s^{-1}(Y\times Z)\cap (U\times O)\subseteq U\times Z$ is
smooth. By Corollary 10.7  in \cite[Ch.\ III]{Ha} the general
fiber of the projection $\cY\to U$ is  smooth as well. In other
words, shrinking $U$ we may assume that all fibers of this
projection are smooth. Since for a point $h=(h_1,\ldots, h_s)\in
U$ the fiber $\cY\cap\pi^{-1}(h)$ maps bijectively via $\Phi_s$
onto $Y\cap (h_1\cdot\ldots\cdot h_s).Z$, (a) follows.

Now (b) follows by an easy induction on $l=\dim Y+\dim Z$, the
case of $l=0$ being trivial. Indeed, applying (a) and the
induction hypothesis to $Y_{\sing}$ and $Z$ and also to $Y$ and
$Z_{\sing}$, for suitable connected algebraic subgroups
$H_1,\ldots, H_s$ and general $(h_1,\ldots ,h_s)\in
H_1\times\ldots\times H_s$ we have that $Y_\reg$ and
$(h_1\cdot\ldots\cdot h_s).Z_\reg$ meet transversally and that
\bals \dim (Y_{\sing}\cap (h_1\cdot\ldots\cdot h_s).Z)&\le
\dim Y_{\sing}+\dim Z-\dim X;\\
\dim (Y\cap (h_1\cdot\ldots\cdot h_s).Z_{\sing})&\le \dim Y+ \dim
Z_{\sing}-\dim X\,. \eals
This immediately implies the desired
result.
\eproof

\subsection{$\G_a$-generated subgroups}\label{specialgr}
%Let $X$ be an algebraic variety.
The following notion is central in the sequel.

\bdefi\label{special} A subgroup $G$ of the automorphism group
$\Aut(X)$ will be called {\em $\G_a$-generated}\/\footnote{
These groups were introduced in the particular case $X=\AA^n$  in
\cite[Definition 2.1]{Po4}, where they were called {\em
$\p$-generated}. Cf.\ also  Definition 1.36 in \cite{Po2} for a
more general notion of an {\em $F$-generated group}.} if it is
generated by a family of one-parameter unipotent subgroups i.e.,
subgroups isomorphic to $\G_a$. \edefi

We give two simple examples.

\bexa\label{fv.12.26.11a} (1) The group $\SAut(X)$  is
$\G_a$-generated. The image of  $\SAut(X)$ under the diagonal
embedding $\SAut(X)\hookrightarrow \SAut(X^m)$ is also a
$\G_a$-generated subgroup.

(2) A connected affine algebraic group acting regularly and
effectively on $X$ is a $\G_a$-generated subgroup of  $\Aut(X)$ if
and only if it does not admit nontrivial characters \cite[Lemma
1.1]{Po2}. In particular, every connected semisimple algebraic group is
$\G_a$-generated. \eexa

It will be important to deal with the infinitesimal
generators of algebraic subgroups of $\Aut(X)$ isomorphic to
$\G_a$. Let us collect the necessary facts.

\bsit\label{remspecial} (1) If the group $\G_a$ acts on an affine
variety $X=\Spec A$ then the associated derivation $\partial$ of
$A$ is locally nilpotent, i.e.\ for every $a\in A$ we can find
$n\in \N$ such that $\partial^n(a)=0$ \cite{Ren}. It is
immediate that for every $f \in \ker \partial$ the derivation $f
\partial $ is again locally nilpotent \cite[\S 1.4, Principle
7]{Fre}.

(2) Conversely, given a locally nilpotent $\kk$-linear derivation
$\partial:A\to A$ and $t\in \kk$, the map $\exp(t\p):A\to A$ is an
automorphism of $A$ \cite[1.1.8]{Fre}.
Furthermore for $\p\neq
0$, $H=\exp(\kk\p)$ is a subgroup of $\Aut(A)$ isomorphic to
$\G_a$. Via the isomorphism
$\Aut(A)\stackrel{\simeq}{\longrightarrow}\Aut(X)$ given by
$g\mapsto (g^{-1})^*$ this yields a one parameter unipotent
subgroup of $\Aut(X)$, which we denote by the same letter $H$.

One can also consider $\p$ as a vector field on $X$. If the ground
field is $\C$ then the action of $H\cong \G_a$ on $X$ is just the
associated phase flow. We often use the term {\em locally
nilpotent vector field} meaning that the associated derivation is
locally nilpotent.

(3) The ring of invariants $\cO (X)^H=\ker\partial$ has
transcendence degree over $\kk$ equal to $\dim X-1$. For any
$H$-invariant function $f\in\cO (X)^H$ the one-parameter unipotent
subgroup $H_f=\exp(\kk f\partial)$ plays an important role in the
sequel (cf.\ \cite{Po6}). It will be called a {\em replica} of
$H$.

(4)  The $H_f$-action has the same general orbits as the
$H$-action. However, the zero locus of $f$ remains pointwise fixed
under the $H_f$-action.
\esit

\bsit\label{no-passaran} In order to illustrate the notions of a
$\G_a$-generated group and a replica, the affine space $X=\AA^n$
is a good choice, as this is done in \cite[\S\S 1,2]{Po4}. For
instance, a replica $H_f$ of a one-parameter unipotent subgroup
$H=\exp(\kk\partial)$ generated by a directional partial
derivative $\partial$ is a one parameter group of shears in the
same direction. For $n=3$, the famous Nagata
automorphism\footnote{Recall that the Nagata automorphism is wild,
see \cite{US}.} is actually a special value of the replica
associated with the locally nilpotent derivation $\p=X\frac{\p}{\p
Y}+Y\frac{\p}{\p Z}$ of the polynomial ring $\kk[X,Y,Z]$ and the
invariant function $f=Y^2-2XZ\in\ker\p$. The problem whether the
subgroup $\SAut(\AA^n)$ coincides with the group of
all automorphisms of $\AA^n$ with Jacobian determinant 1 is
still widely open (see \cite[Problem 2.1, Examples 2.3 and
2.5]{Po4}; cf.\ also \cite[Proposition 9]{FurLa}). Recall  that this is indeed the case in dimension 2 due to the Jung-van der Kulk Theorem, see {\em ibid}.
 \esit

Given an algebraic variety we denote by $\LND(X)$ the set of all
locally nilpotent vector fields on $X$. If $G\subseteq \Aut(X)$ is
any subgroup then the vector fields in $\LND(X)$ generating
one-parameter unipotent subgroups of $G$ form a subset $\LND(G)$
of $\LND(X)$. This set is a cone (i.e., $\kk\cdot\LND(G)\subseteq
\LND(G)$) stable under conjugations $\p \mapsto g^* \p
(g^*)^{-1}$, where the automorphism $g^*: A\to A$ is induced by $g
\in G$ (cf. \cite[\S 1.4, Principle 1d]{Fre}).

Let now $G\subseteq \SAut(X)$ be a $\G_a$-generated subgroup.
In the sequel we consider subsets $\cN\subseteq \LND(G)$ of
locally nilpotent vector fields such that the associated
one-parameter subgroups $\cG=\{\exp(\kk\p):\p\in \cN\}$ form a
generating set of algebraic subgroups for $G$. By abuse of
language, we often say that $\cN$ is a generating set of $G$, and
we write $G=\langle\cN\rangle$.

>From Proposition \ref{1.5} we deduce the following result.

\bcor\label{locclo} Given a $\G_a$-generated subgroup
$G=\langle\cN\rangle$ of $\Aut(X)$, where $\cN\subseteq\LND(G)$ is
stable under conjugation in $G$, there are locally nilpotent
vector fields $\p_1,\ldots,\p_s\in\cN$ which span the tangent
space $T_p (G.p)$ at every point $p\in X$. \ecor

For a point $p\in X$ we let $\LND_p(G)\subseteq T_p X$ denote the
{\em nilpotent cone} of all tangent vectors $\p(p)$, where $\p$
runs over $\LND(G)$. By Corollary \ref{locclo} we have $T_p
(G.p)=\Span \LND_p(G)$.\footnote{Cf.\ Corollary \ref{5.20} below.}
Thus a point $p\in X_{\reg}$ is $G$-flexible (see Definition
\ref{gflex}) if and only if the cone $\LND_p(G)$ spans the whole
tangent space $T_p X$ at $p$.

Applying Corollaries \ref{gflex1} and \ref{locclo} to the
special automorphism group $G=\SAut(X)$ yields the equivalence
(i)$\Leftrightarrow$(iii) in Theorem \ref{mthm} in the
Introduction.

\bcor\label{tr-fl} Given an affine variety $X$
the action of $\SAut
(X)$ on $X_{\rm reg}$ is transitive if and only if $X$ is
flexible.
\ecor

\section{Infinite transitivity}\label{inftrans}

\subsection{Main theorem}\label{maint-m}
In this section we show that the special automorphism group of a
flexible irreducible affine variety $X$ acts infinitely
transitively on $X_\reg$. We state this in a more general setup
which turns out to be necessary for later applications. Let us
first introduce the following useful notation.

\bdefi\label{2.1}  Let $X$
be an irreducible affine algebraic variety. A set $\cN$ of
locally nilpotent vector fields on $X$ is said to be {\em
saturated} if it satisfies the following two conditions.

\bnum
\item $\cN$ is closed under conjugation by elements in $G$, where $G$
is the subgroup of $\SAut(X)$ generated by $\cN$.
\item $\cN$ is closed under taking replicas,
i.e.\ for all $\p\in \cN$ and $f\in \ker \p$ we have $f\p\in \cN$.
\enum
\edefi

We note that replicas appear implicitely at many places of the literature, see for instance \cite{KZ1};  see also \cite[Definition 2.1]{Po4} for a related definition.

Clearly every collection $\cN^o$ of locally nilpotent
vector fields on $X$ can be extended to a saturated set $\cN$. We
note that in general this extended set generates a much larger
group than $\cN^o$. For instance, if $X=\AA^2$ is equipped with
coordinates $(x,y)$ then  $\cN^o=\{\partial/\partial y\}$
generates the translations $(x,y)\mapsto (x, y+t)$ while its
saturation $\cN$ generates the group of all shears $(x,y)\mapsto
(x,y+P(x))$ with $P\in \kk[x]$.

The next result implies Theorem \ref{mthm} in the Introduction.

\bthm\label{maincor} Let $X$ be an irreducible affine
algebraic variety of dimension $\ge 2$ and let  $G\subseteq
\SAut(X)$ be a subgroup generated by a saturated set $\cN$ of
locally nilpotent vector fields, which acts with an open orbit
$O\subseteq X$. Then $G$ acts infinitely transitively on $O$.
\ethm

Before starting the proof let us mention the following
interesting class of examples.

 \bexa\label{gisu} 1. ({\em Gizatullin surfaces.})
 These are normal
affine surfaces which admit a completion by a chain of smooth
rational curves. Due to Gizatullin's Theorem (\cite[II, Theorems 2
and 3]{Gi}\footnote{In \cite[II]{Gi} the result is stated in terms
of $\Aut (X)$, but the proof applies to $\SAut (X)$.}; see also
\cite{Du}) a normal affine surface $X$ different from $\AA^1\times
(\AA^1\setminus\{0\})$ is Gizatullin if and only if the special
automorphism group $\SAut(X)$ has an open orbit with a finite
complement. By Theorem \ref{maincor} the group $\SAut (X)$ acts
infinitely transitively on this orbit.
It was conjectured in \cite[II]{Gi} that in zero characteristic
this orbit coincides with $X_{\reg}$ i.e., that every Gizatullin
surface is flexible. This is definitely not true in positive
characteristic, where the automorphism group $\Aut(X)$ of a
Gizatullin surface $X$ can have fixed points that are regular
points of $X$ \cite{DG}. The Gizatullin Conjecture is true for the
Gizatullin surfaces given in $\AA^3$ by equations\footnote{Over an
arbitrary base field of characteristic zero.} $xy-f(z)=0$; see
\cite[Theorem 3.1]{AKZ} and \cite{ML1, ML2}. Yet another class
of flexible Gizatullin surfaces consists of the Danilov-Gizatullin
surfaces, see \cite{Gi0} (see also \cite[Theorem 5]{Do1}).
We refer the reader to \cite{FKZ} and the references therein
for a study of one-parameter groups acting on Gizatullin surfaces.
\eexa

\bsit\label{nz}
For  subsets $\cN\subseteq \LND(X)$ and  $Z\subseteq X$
we let $\cN_Z=\{ \p\in\cN: \p|Z=0\}$ be the set of locally
nilpotent vector fields in $\cN$ vanishing on $Z$.
If $G=\langle\cN\rangle$ is the group generated by $\cN$ then
$\cN_Z$ generates a subgroup denoted
\be\label{GNZ}
G_{\cN,Z}=\langle H=\exp(\kk\p): \p\in \cN_Z \rangle \subseteq G\,.\ee
Clearly the automorphisms in $G_{\cN,Z}$ fix the set $Z$
pointwise. In the case $\cN=\LND(G)$ we simply write $G_Z$ instead
of $G_{\cN,Z}$.

We emphasize that if $\cN$  is saturated then so is $\cN_Z$.
Hence in this case  the group
$G_{\cN,Z}=\langle\cN_Z\rangle$ is again generated by a saturated
set of locally nilpotent derivations.\esit

With these notations our main technical result can
be formulated as follows.

\bthm\label{mainnew} Let $X$ be an irreducible affine
algebraic variety of dimension $\ge 2$ and let $G\subseteq
\SAut(X)$ be a subgroup generated by a saturated set $\cN$ of
locally nilpotent vector fields, which acts with an open orbit
$O\subseteq X$. Then for every finite subset $Z\subseteq O$ the
group $G_{\cN,Z}$ acts transitively on $O\backslash Z$. \ethm

Before embarking on the proof let us show how Theorem
\ref{maincor} follows from Theorem \ref{mainnew}.

\bproof[Proof of Theorem \ref{maincor}] Let $x_1,\ldots ,x_m$ and
$x_1',\ldots, x_m'$ be sequences of points in $O$ with $x_i\ne
x_j$ and $x_i'\ne x_j'$ for $i\ne j$. Let us show by induction on
$m$ that there is an automorphism $g\in G$ with $g.x_i=x_i'$ for
all $i=1,\ldots ,m$. As $G$ acts transitively on $O$ this is
certainly true for $m=1$. For the induction step suppose that
there is already an automorphism $\alpha\in G$ with $\alpha
.x_i=x_i'$ for $i=1,\ldots, m-1$. Applying Theorem \ref{mainnew}
to $Z=\{x_1',\ldots, x_{m-1}'\}$ we can also find an
automorphism $\beta\in G_{\cN,Z}$ with $\beta
(\alpha(x_m))=x_m'$.  Clearly then $g=\beta\circ\alpha$
satisfies $g.x_i=x_i'$ for all $i=1,\ldots ,m$.
\eproof

\subsection{Proof of Theorem \ref{mainnew}}
\label{prmthm} To deduce Theorem \ref{mainnew} we need a few
preparations. As before $X$ stands for an irreducible affine
algebraic variety.  Let us introduce the following technical
notion.

\bdefi\label{orbitsep} Let $G\subseteq\SAut(X)$ be a
$\G_a$-generated subgroup and let $\Omega\subseteq X$ be a subset
invariant under the $G$-action, i.e.\ $G.\Omega\subseteq \Omega$.
We say that a locally nilpotent vector field $\p\in \LND(G)$ with
associated one-parameter subgroup $H=\exp(\kk\p)$ satisfies the
{\em orbit separation property} on $\Omega $,  if there is an
$H$-stable subset $U(H)\subseteq \Omega$ such that \bnum[(a)]
\item for each $G$-orbit $O$
contained in $\Omega$, the intersection $U(H)\cap O$ is open and
dense in $O$, and

\item the global $H$-invariants $\cO(X)^H$
separate all one-dimensional $H$-orbits in $U(H)$. \enum \noindent

The reader should note that we allow $U(H)\cap O$ to contain or
even to consist of $0$-dimensional $H$-orbits. We also emphasize
that $\Omega$ can be {\em any} union of orbits and can e.g.\
contain orbits in the singular part of $X$. As a trivial case, if $\p=0$ then $H=\{1\}$ and the orbit separation property is trivially satisfied with $U(H)=\Omega$.

Similarly we say that a set of locally nilpotent vector fields
$\cN$ satisfies the {\em orbit separation property} on $\Omega$ if
this holds for every $\p\in \cN$. \edefi

As we shall see in Example \ref{conter} the orbit separation
property is not necessarily satisfied on every $G$-stable subset.
However, there are interesting geometric situations where this
property holds for arbitrary subsets stabilized by $G$, see
Subsection \ref{matvar}. In the following remarks we indicate
possible choices of a good set $\Omega$.

\brems\label{exanew} 1. Let $\p$ be a locally nilpotent vector field on $X$ and let $H=\exp(\kk \p)$ be the subgroup of $\SAut(X)$ generated by $\p$. According to  \cite[Theorem 3.3]{PV} the field of rational invariants $\mathcal{M}er(X)^H$ is the quotient field of the ring  $\cO(X)^H$ of regular invariants. Hence by a corollary of the Rosenlicht theorem on rational invariants (see \cite[Proposition 3.4]{PV})  the regular invariants $\cO(X)^H$ separate orbits on an $H$-invariant open dense subset $U(H)$ of $X$. With such a set $U(H)$ condition (b) in Definition \ref{orbitsep} is automatically satisfied.

If furthermore $\Omega=O$ is an open $G$-orbit in $X$ then $U(H)$ can be chosen to be contained in $O$ so that the other requirements of Definition  \ref{orbitsep} are as well satisfied. It follows that {\em every $\p\in \LND(G)$ satisfies the orbit separation property on an open $G$-orbit $\Omega=O$ in $X$.}

2. In a similar fashion, given a  locally nilpotent derivation
$\p\in\LND(G)$, it satisfies the orbit separation property on any
set $\Omega$ which is a union of $G$-orbits meeting $U(H)$, where
$H=\exp(\kk\p)$. In particular this property holds for general
$G$-orbits (cf.\ Corollary \ref{1.4}).

3. Suppose that $\Omega$ in Definition \ref{orbitsep} consists of a single $G$-orbit $O$. Let $\p\in \LND(G)$ and $U(H)$ be  as in  \ref{orbitsep}. Shrinking $U(H)$
if necessary we can achieve that $U(H)$ has a geometric quotient $U(H)/H$, which admits a locally closed embedding into some $\AA^N$
by regular invariants in $\cO(X)^H$ (cf.\ also \cite[Theorem 4.4]{PV}).
\erems

We need the following simple Lemma.

\blem \label{3.2} If a locally nilpotent vector field $\p\in
\LND(G)$ satisfies the orbit separation property on a $G$-stable
subset $\Omega\subseteq X$ then also every replica $f\p$, $f\in
\ker\p$, and every $g$-conjugate $g^*(\p)=g\circ\p \circ g^{-1}$,
$g\in G$, has this property. \elem

\bproof Let $\p'=f\p$ be a replica of $\p$ with associated
one-parameter subgroup $H'$. In the case $f=0$ the assertion is
obvious. Otherwise the one dimensional orbits of $H'$ are also one
dimensional orbits of $H$, and the $H$ and $H'$ invariant
functions are the same. Hence setting $U(H')=U(H)$, (a) and (b) in
Definition \ref{orbitsep} are again satisfied for $H'$. The fact
that any $g$-conjugate of $\p$ has again the orbit separation
property can be left to the reader. \eproof

For the remaining
part of this subsection we fix the following notation.

\bsit\label{not2.2}
Let $G\subseteq\SAut(X)$ be a $\G_a$-generated
subgroup generated by a saturated set $\cN$ of locally nilpotent
vector fields. Let $\Omega\subseteq X$ be a $G$-stable subset. We
choose $\p_1,\ldots, \p_s\in \cN$ with associated one-parameter
subgroups $H_\sigma=\exp(\kk\p_\sigma)$. We assume in Lemma
\ref{zaru} below that the following two conditions are
satisfied. \bnum[(1)]
\item $\p_1,\ldots, \p_s\in \cN$ span $T_x(G.x)$
for every point $x\in \Omega$ (see \ref{locclo}), and
\item $\p_\sigma$
has the orbit separation property on $\Omega$ for all
$\sigma=1,\ldots, s$. \enum Consequently there are subsets
$U(H_\sigma)\subseteq \Omega$ such that conditions (a) and (b)  in
Definition \ref{orbitsep} are satisfied with $H=H_\sigma$. We let
$$
V=\bigcap_{\sigma=1}^s U(H_\sigma)\,.
$$
In particular,
\bnum[(i)]
\item$V\cap O$ is open and dense in $O$
for every orbit $O$ contained in $\Omega$, and
\item any two points in $V$ in different one dimensional
$H_\sigma$-orbits can be separated by an $H_\sigma$-invariant
function on $X$ for all $\sigma=1,\ldots, s$.
\enum
\esit

\blem\label{zaru}  With the notation and assumptions as in
\ref{not2.2} above, for any pair of distinct points  $x,y\in
\Omega$ lying in $G$-orbits of dimension $\ge 2$ there exists
an automorphism $g\in G$ such that

(a) $g.x,\ g.y \in V$, and

(b) $g.x$ and $g.y$ are lying in different
$H_\sigma$-orbits\footnote{Possibly of dimension 0.} for
all $\sigma=1,\ldots, s$.
\elem

\bproof
(a) Since $G$ acts transitively on every $G$-orbit
$O$ in $\Omega$ and $V\cap O$ is dense in $O$,
we can find $g\in G$ with $g.x\in V$.
Replacing $x$ by $g.x$ we may assume that $x\in V$.
For some $\sigma\in\{1,\ldots, s\}$
we have $H_\sigma.y\cap V\ne\emptyset$.
Taking $h\in H_\sigma$ general we have $h.
x\in V$ and $h.y\in V$, as required.

(b) By (a) we may assume that $x,y\in V$. The property that  $g.x$
and $g.y$ are in different $H_\sigma$-orbits is an open
condition for $(x,y)$ running over the space $V\times V$.
Thus by recursion it suffices to find $g\in G$ such that (b) is
satisfied for a fixed $\sigma$. If $x$ and $y$ are already in
different $H_\sigma$-orbits then there is nothing to show.

So suppose that this is not the case
and so $x,y$ are sitting on the same $H_\sigma$-orbit,
which is then necessarily one dimensional.
By assumption the vector fields $\p_1,\ldots, \p_s$
span the tangent space $T_{x}(G.x)$ at $x$,
and the $G$-orbit of $x$ has dimension $\ge 2$.
Hence  $\p_\tau$ is not tangent to $H_\sigma.x$ at $x$
for some $\tau$. In particular the orbits $H_\sigma.x$
and $H_\tau.x$ are both of dimension one and
have only finitely many points in common.

If $x$ and $y$ are in different $H_\tau$-orbits
then we can choose a global $H_\tau$-invariant $f$ with $f(x)=1$
and $f(y)=0$. The group $H=\exp(\kk f\p_\tau)$
is contained in $G$,  fixes $y$ and moves $x$
along $H_\tau.x$. Hence for a general $g\in H_\tau $
the points $g.x$ and $g.y=y$ lie on different $H_\sigma$-orbits.

Assume now that $x$ and $y$ belong to the same $H_\tau$-orbit. We
claim that again $g.x$ and $g.y$ are in different
$H_\sigma$-orbits for a general $g\in H_\tau$.

To show this claim we consider  $h_t=\exp(t\p_\tau)\in H_\tau$. By
assumption $h_a.x=y$ for some $a\in \kk$, $a\neq 0$. We can find
an $H_\sigma$-invariant function $f$ on $X$, which induces a
polynomial $p(t)=f(h_t.x)$ of positive degree in $t\in \kk$. If
$h_t.x$ and $h_t.y$ are in the same $H_\sigma$-orbits for a
general $t\in \kk$ then
$$
p(t)=f(h_t.x)=f(h_t.y)=f(h_t.(h_a.x))=f(h_{a+t}.x)=p(t+a),
$$
which is impossible. Hence for a general $g=h_t\in H_\tau$
the points $g.x$ and $g.y$ are in different $H_\sigma$-orbits,
as desired.
\eproof

\blem\label{new1}
With the notations as in \ref{not2.2}
assume that $x,y\in V$ are distinct points lying in different
(possibly zero dimensional)
$H_\sigma$-orbits for all $\sigma=1,\ldots, s$.
Then the vector fields $\p\in \cN$ vanishing
at $x$ span $T_y(G.y)$.
\elem

\bproof
The vectors $\p_\sigma(y)$ with $1\le \sigma\le s$  span
the tangent space $T_y(G.y)$. Thus it suffices to
find replicas $\p_1', \ldots, \p_s'$
of $\p_1,\ldots , \p_s$, which vanish at $x$
and are equal to $\p_\sigma$ at the point $y$.

If the $H_\sigma$-orbit of $x$ is a point,
then necessarily $\p_\sigma$ vanishes at $x$ and
we can choose $\p_\sigma'=\p_\sigma$. If the $H_\sigma$-orbit
of $y$ is a point then $\p_\sigma(y)=0$ and so
we can take $\p_\sigma'=0$. Assume now that both
$H_\sigma.x$ and $H_\sigma.y$ are one dimensional.
By our construction of $V$ there is an $H_\sigma$-invariant
function $f_\sigma$ on $X$ with $f_\sigma(x)=0$ and $f_\sigma(y)=1$.
Hence $\p_\sigma'=f_\sigma\p_\sigma$ is a
locally nilpotent vector field on $X$ vanishing in $x$
and equal to $\p_\sigma$ at $y$.
\eproof

\bcor\label{2.8} For each $x\in \Omega$ and every $G$-orbit
$O\subseteq \Omega$ the group $G_{\cN,x}$ as in \ref{nz}
acts transitively on $O\backslash\{x\}$.\footnote{In particular,
it is transitive on $O$ if $x\not\in O$. } \ecor

\bproof Let $y$ be a point in $O\backslash\{x\}$. With the
notations as in \ref{not2.2}, according to Lemma \ref{zaru} there
is an automorphism $g\in G$ with $g.x,\ g.y\in V$ such that $g.x$,
$g.y$ are in different (possibly $0$-dimensional)
$H_\sigma$-orbits for $i=1,\ldots, s$. By Lemma \ref{new1} the
vector fields $\p\in \cN$ vanishing at $g.x$ span $T_{g.y}(O)$.
Using the fact that $\cN$ is stable under conjugation by elements
$g\in G$ it follows that the vector fields in $\cN$ vanishing at
$x$ span the tangent space $T_y(O)$.  In other words, $y$ is a
$G_{\cN,x}$-flexible point on the orbit closure $\bar O$. Applying
\ref{gflex1} we obtain that $G_{\cN,x}$ acts transitively on
$O\backslash\{x\}$. \eproof

\bproof[Proof of Theorem \ref{mainnew}] By Remark \ref{exanew}(1)
the orbit separation property is satisfied on the open orbit
$\Omega=O$. Given a set $Z=\{x_1,\ldots, x_m\}\subseteq O$
 of $m$ distinct points we consider $Z_\mu=\{x_1,\ldots,x_\mu\}$
for $\mu=1,\ldots, m$. Let us show by induction on $\mu$ that
$G_{\cN,Z_\mu}$ acts transitively on $O\backslash Z_\mu$. For
$\mu=1$ this follows from Corollary \ref{2.8}. Assuming for some
$\mu <m$ that $G_{\cN,Z_\mu}$ acts transitively on $O\backslash
Z_\mu$,  Corollary \ref{2.8} implies that
$(G_{\cN,Z_\mu})_{x_{\mu +1}}=G_{\cN,Z_{\mu +1}}$ acts
transitively on $O\backslash Z_{\mu+1}$. \eproof

\subsection{Examples of non-separation of orbits}
Suppose as before that a subgroup $G\subseteq\SAut (X)$ is
generated by a saturated set $\cN$ of locally nilpotent vector
fields. Then $\cN$ satisfies the orbit separation property
\ref{orbitsep} on a general $G$-orbit (see Remark
\ref{exanew}(2)), while this is not always true on an arbitrary
$G$-orbit. Furthermore, the following example shows that on the
union of two $G$-orbits this property might fail although it is
satisfied on every single orbit.

\bexa\label{non-balanced} On the affine 4-space
$\AA^4=\Spec\kk[X,Y,Z,U]$  let us consider the locally nilpotent
vector fields
$$\p_1=Y\frac{\p}{\p X}+Z\frac{\p}{\p Y}\quad\mbox{and}\quad
\p_2=\frac{\p}{\p U}\,.$$ Let $G\subseteq\SAut(\AA^4)$ be the
$\G_a$-generated subgroup generated by $\p_1,\p_2$ and all their
replicas, and let $\cN\subseteq \LND(G)$ denote the saturated set
generated by $\p_1$ and $\p_2$.  Note that $O_\pm=\{Y=\pm
1,\,Z=0\}$ is naturally isomorphic to the $(X,U)$-plane with the
restrictions of $\p_1$ and $\p_2$ given by $\pm \frac{\p}{\p X}$
and $\frac{\p}{\p U}$, respectively. Therefore $O_\pm$ is a
$G$-orbit. It is easily seen that $\ker\p_1=\kk[Z,\,Y^2-2XZ,U]$
and $\ker \p_2=\kk[X,Y,Z]$. Hence the $G$-orbits $O_+$ and
$O_-$  are not separated by $H_i$-invariants, where
$H_i=\exp(\kk\p_i)$, $i=1,2$. In particular, $\cN$ does not
satisfy the orbit separation property on $O_+\cup O_-$. However,
this property is satisfied on $\Omega=O_+$ and also on
$\Omega=O_-$ separately as this is the case for $\p_1$ and $\p_2$
(cf.\ Lemma \ref{3.2}).

We note also that the isomorphism $\sigma:O_+\to O_-$ with
$\sigma(x,1,0,u)=(-x,-1,0,u)$ commutes with the actions of $H_1$
and $H_2$. Hence there is no collective transitivity on $O_+\cup
O_-$ in the sense of Theorem \ref{maincorstrat} below, while $G$
acts on every single orbit $O_\pm$ indeed infinitely transitively.
\eexa

According to our next example one cannot expect infinite
transitivity of $G$ on an arbitrary $G$-orbit $O$ of dimension
$\ge 2$ without assuming the orbit separation property on $O$.
However, compare Theorem \ref{maincorstrat} below for a
positive result.

\bexa\label{conter}
Consider the locally nilpotent
derivations
$$\p_1=Y\frac{\p}{\p X}+Z\frac{\p}{\p Y}+U\frac{\p}{\p
Z}\quad\mbox{and}\quad \p_2=Z\frac{\p}{\p X}+X\frac{\p}{\p
Y}+U\frac{\p}{\p Z}\,$$
of the polynomial ring $\kk[X,Y,Z,U]$.
We claim that
$\ker\p_1=\kk[p_1,p_2,p_3,p_4]$, where
\bals
&p_1=U,\quad p_2=Z^2-2YU,\quad p_3=Z^3-3YZU+3XU^2\,,\quad \mbox{and}\\
&
p_4=\frac{p_2^3-p_3^2}{p_1^2}=9X^2U^2-18XYZU+6XZ^3-3Y^2Z^2+8Y^3U\,.
\eals
Indeed, the image of the map
$$\rho=(p_1,\ldots,p_4):\AA^4\to \AA^4\,$$
 is contained in
 the hypersurface $$F=\{X_1^2X_4-X_2^3+X_3^2=0\}$$
which is singular along the line
$F_{\rm sing}=\{X_1=X_2=X_3=0\}$. Being regular in
codimension one, $F$ is normal.
We have
$$\bar 0\in F_{\rm sing}\quad\mbox{and}\quad
\rho^{-1}(\bar 0)=\rho^{-1}(F_{\rm sing})=\{Z=U=0\}
=:L\subseteq\AA^4\,.$$ By the Weitzenb\"ock Theorem (see e.g.
\cite{Kr})
%there exists a categorical quotient
$E=\Spec (\ker\p_1) $
is an affine algebraic variety. The
inclusions
$$\kk[p_1,p_2,p_3,p_4]\subseteq\ker\p_1
\subseteq\kk[X,Y,Z,T]$$
lead to morphisms
$$\AA^4\stackrel{\pi}{\longrightarrow} E
\stackrel{\mu}{\longrightarrow}  F\,, \qquad\mbox{where} \quad
\mu\circ\pi=\rho\,.$$ We claim that  $\mu$ is an  isomorphism.
Since both $E$ and $F$ are normal affine threefolds, by the
Hartogs Principle \cite[Proposition 7.1]{Da} $\mu$ is an
isomorphism if it is so in codimension one. In turn, it suffices
to check that $\mu$ admits an inverse morphism defined on
$F_{\reg}$. The latter follows once we know that $\rho$ separates
the $H_1$-orbits in $\AA^4$ outside the plane $L$  and that
$F_{\reg}=\rho(\AA^4\backslash L)$. Indeed, then $\pi$ separates
them as well, and $\mu$ induces a bijection between
$\pi(\AA^4\backslash L)\subseteq E$ and $F_{\reg}$.

%According to our convention in \ref{remspecial}(2) t
The action of $H_1$ on $\AA^4$ is given by
\be\label{haction} (-t).\left(%
\begin{array}{c}
  x  \\
  y  \\
  z  \\
  u \\
\end{array}%
\right)=
\left(%
\begin{array}{c}
  x + ty +\frac{t^2}{2}z+\frac{t^3}{6}u\\
  y +tz+\frac{t^2}{2}u \\
  z +tu \\
  u \\
\end{array}%
\right)\,.\ee
Let $O$ be an $H_1$-orbit contained in
$\AA^4\backslash L$. Suppose first that $p_1\vert O=U\vert O=u\neq
0$. Letting $t=-z/u$ in (\ref{haction}) we get a point
$A=(x,y,0,u)\in O$. Since
$$p_2(A)=-2yu\quad\mbox{and}\quad p_3(A)=3xu^2$$
we can recover the coordinates \be\label{co1} y=-(p_2\vert
O)/2u\quad\mbox{and}\quad x=(p_3\vert O)/3u^2\,.
\ee
Thus $O$ is uniquely determined by the image $\rho(O)\in F$.

 Suppose further that $p_1\vert O=U\vert O=0$.
 Since $O\cap L=\emptyset$
 then $Z\vert
 O=z\neq 0$. Taking $t=-y/z$ in
(\ref{haction}) yields a point $A=(x,0,z,0)\in O$.
 Since
 $$p_2(A)=z^2,\,p_3(A)=z^3,\quad\mbox{and}\quad  p_4(A)=6xz^3 $$
  we can recover the values\footnote{Formulas (\ref{co1})
 and (\ref{co2}) define sections of
 $\rho$ in the open sets $U\neq 0$ and $Z\neq 0$, respectively.
 This shows
 that $\rho: \AA^4\backslash L\to F_{\reg}$ is a principal
 $\AA^1$-bundle.}
 \be\label{co2}
 z=(p_3\vert O)/(p_2\vert O)\quad\mbox{and}\quad x=
 (p_4\vert O)/6z^3\,. \ee
Now both claims
follow.

The generators $p_1,\ldots,p_4$ of the algebra of $H_1$-invariants
vanish on the plane $L=\{Z=U=0\}$ so that every $H_1$-invariant is
constant on $L$.  Since $\p_2$ is obtained from $\p_1$ by
interchanging $X,Y$, by symmetry also every $H_2$-invariant is
constant on $L$. Letting $G=\langle\Sat(H_1,H_2)\rangle$ be the
subgroup generated by $H_1$, $H_2$, and all their replicas,  it
easily follows that $G$ stabilizes $L$ and that the $G$-action on
$L$ factors through the $\SL_2$-action   associated with the
$\sgoth \lgoth_2$-algebra generated by the vector  fields $\p_1|_L
=Y\frac{\p}{\p X} $ and $\p_2|_{L}=X\frac{\p}{\p Y}$. In
particular, the action of $G$ on its orbit $L\backslash\{0\}$ is
not even $2$-transitive, the linear dependence being an
obstruction.

Observe finally that the three dimensional $G$-orbits in $\AA^4$
are separated by the $G$-invariant function $U$. \eexa

\subsection{$G_Y$-orbits}

Let as before $G\subseteq \SAut(X)$ be generated by a saturated
set $\cN$ of locally nilpotent vector fields on an
irreducible affine variety $X=\Spec A$. For a subvariety
$Y\subseteq X$ we consider the subgroup $G_{\cN,Y}$ of $G$ as
defined in (\ref{GNZ}).

\bthm\label{2.9}
Under notation as above, assume that  $\dim
X\ge 2$ and that $G$ acts on $X$ with an orbit $O$ open in $X$.
Then letting $Y= X\backslash O$ the subgroup $G_{\cN,Y}$ acts
transitively and hence infinitely transitively on $O$. \ethm

\bproof Since $G_{\cN,Y}$ is generated by a saturated set of
locally nilpotent derivations, by Theorem \ref{maincor} it
suffices to show that $G_{\cN,Y}$ acts transitively on $O$.

Using Corollary \ref{locclo} we can choose $\p_1,\ldots,
\p_s\in\cN$ spanning the tangent space $T_xX$ at each point $x\in
O$. Letting $I$ denote the ideal of $Y$ in $A$, we claim that for
every $\sigma=1,\ldots,s$ there is a nonzero function $f_\sigma\in
I\cap\ker\p_\sigma$. Let $H_\sigma=\exp(\kk\p_\sigma)\subseteq G$,
$\sigma=1,\ldots,s$. The set $Y$ being $G$-invariant, for every
nonzero function $f\in I$ the orbit $H_\sigma.f$ spans in $A$ an
$H_\sigma$-invariant finite dimensional subspace $E_\sigma$
contained in $I$. By the Lie-Kolchin Theorem there is a nonzero
element $f_\sigma\in E_\sigma$ fixed by $H_\sigma$. This proves the
claim.

Let $p\in X$ be a general point so that $f_\sigma(p)\ne 0$ for
$\sigma=1,\ldots, s$. We can normalize the invariants $f_\sigma$
so that $f_\sigma(p)=1$ and $f_\sigma|Y=0$. The derivation
$f_\sigma\p_\sigma$ then vanishes on $Y$ and so the replica
$H_{\sigma,f_\sigma}=\exp(\kk f_\sigma\p_\sigma)\subseteq
G_{\cN,Y}$
 of $H_\sigma$
fixes $Y$ pointwise while moving $p$ in the direction of
$\p_\sigma(p)$. It follows that the $G_{\cN,Y}$-orbit of $p$ is
open in $O$.

Let now $q\in O$ be an arbitrary point. Choose $g\in G$ with
$g.p=q$. Since $g$ stabilizes $Y$ the subgroup
$H'_{\sigma,f_\sigma}=gH_{\sigma,f_\sigma} g^{-1}\subseteq
G_{\cN,Y}$ fixes $Y$ pointwise and moves $q$ into the direction of
${d}g(\p_\sigma)(q)$. It follows that also the $G_{\cN,Y}$-orbit
of $q$ is open. Finally $G_{\cN,Y}$ has $O$ as an open
orbit.
\eproof

\brem\label{Gire}\footnote{ We are grateful to
M.H.~Gizatullin for this observation. } The conclusion of the
theorem remains valid if we replace the group $G_{\cN,Y}$ by its
subgroup $G^m_{\cN,Y}$ generated by the replicas in $\cN$ which
vanish on the $m$th infinitesimal neighborhood of $Y$ in $X$,
where $m\in\N$. Indeed, it suffices to replace the functions
$f_\sigma$ in the proof by their $m$th powers
$f^m_\sigma$.\erem

\section{Collective infinite transitivity}\label{colinftr}

\subsection{Collective transitivity on $G$-varieties}\label{gvarty}
By {\em collective infinite transitivity} we mean a possibility to
transform simultaneously (that is, by the same automorphism) an
arbitrary finite set of points along their orbits into some given
position. Applying the methods developed in Section 2 we can
deduce the following generalization of Theorem \ref{mainnew}.
Below $X$ stands as usual for an irreducible affine algebraic
variety.

\bthm\label{maincorstrat} Let $G\subseteq \SAut(X)$ be a subgroup
generated by a saturated set $\cN$ of locally nilpotent vector
fields, which has the orbit separation property on a $G$-invariant
subset $\Omega\subseteq X$. Suppose that $x_1,\ldots ,x_m$ and
$x_1',\ldots, x_m'$ are points in $\Omega$ with $x_i\ne x_j$ and
$x_i'\ne x_j'$ for $i\ne j$ such that for each $j$ the orbits
$G.x_j$ and $G.x_j'$ are equal and of dimension $\ge 2$. Then
there exists an element $g\in G$ such that $g.x_j=x'_j$ for
$j=1,\ldots, m$. \ethm

As in Section 2 this will be deduced from the following
more technical result.

\bthm \label{mainnewstrat} Let $G\subseteq \SAut(X)$ be a subgroup
generated by a saturated set $\cN$ of locally nilpotent vector
fields. Suppose that $\cN$ has the orbit separation property on a
$G$-invariant subset $\Omega$. If $Z\subseteq \Omega$ is a finite
subset and $O\subseteq \Omega$ is an orbit of dimension $\ge 2$,
then the group $G_{\cN,Z}$ acts transitively on $O\backslash Z$.
\ethm

\bproof With $Z_\mu=\{x_1,\ldots,x_\mu\}$ let us show by induction
on $\mu$ that $G_{\cN,Z_\mu}$ acts transitively on $O\backslash
Z_\mu$ for every $G$-orbit $O\subseteq \Omega$ of dimension $\ge
2$. For $\mu=1$ this is just Corollary \ref{2.8}. Assuming for
some $\mu <m$ that $G_{\cN,Z_\mu}$ acts transitively on
$O\backslash Z_\mu$, Corollary \ref{2.8} also implies that
$(G_{\cN,Z_\mu})_{x_\mu +1}=G_{\cN,Z_{\mu +1}}$ acts transitively
on $O\backslash Z_{\mu+1}$. Note that by Lemma \ref{3.2} at each
step the set $\cN_{Z_\mu}$ has again the orbit separation property
on $\Omega$ so that Corollary \ref{2.8} is indeed applicable.
\eproof

\bproof[Proof of Theorem \ref{maincorstrat}] As in the proof
Theorem \ref{maincor} we proceed by induction on $m$, the case
$m=1$ being trivial. For the induction step suppose that there is
already an automorphism $\alpha\in G$ with $\alpha .x_i=x_i'$ for
$i=1,\ldots, m-1$. Applying Theorem \ref{mainnewstrat} to
$Z=\{x_1',\ldots, x_{m-1}'\}$ we can also find an automorphism
$\beta\in G_{\cN,Z}$ with $\beta (\alpha(x_m))=x_m'$.
Clearly then $g=\beta\circ\alpha$ has the required property.
\eproof

\subsection{Infinite transitivity on matrix varieties}\label{matvar}
In this subsection we apply our methods in a concrete setting
where $X=\Mat(n,m)$ is the set of all $n\times m$ matrices over
$\kk$ endowed with the natural stratification by rank. We
assume below that $mn\ge 2$. Let us introduce the following
terminology.

Let $X_r\subseteq X$ denote the subset of matrices of rank $\le
r$. The product $\SL_n\times\SL_m$ acts naturally on $X$ via the
left-right multiplication preserving the strata $X_r$. For every
$k\neq l$ we let $E_{kl}\in\sl_n$ and $E^{kl}\in\sl_m$ denote the
nilpotent matrices with $x_{kl}=1$ and the other entries equal
zero\footnote{Notice that $E_{kl}=E^{kl}$ if $m=n$.}. Let further
$H_{kl}=I_n+\kk E_{kl}\subseteq \SL_n$ and $H^{kl}=I_m+\kk
E^{kl}\subseteq \SL_m$ be the corresponding one-parameter
unipotent subgroups in the first and the second factor of
$\SL_n\times\SL_m$, respectively, acting on the stratification
$X=\bigcup_r (X_r\setminus X_{r-1})$ in a natural way. We also let
$\delta_{kl}$ and $\delta^{kl}$, respectively, denote the
corresponding locally nilpotent vector fields on $X$ tangent to
the strata.

We call {\em elementary} the one-parameter unipotent subgroups
$H_{kl}$, $H^{kl}$, and all their replicas. In the following
theorem we establish the collective infinite transitivity on the
above stratification of the subgroup $G$ of $\SAut(X)$ generated
by the two sides elementary subgroups (cf.\ \cite{Re1}).

By a well known theorem of linear algebra, the subgroup
$\SL_n\times\SL_m\subseteq G$ acts transitively on each stratum
$X_r\setminus X_{r-1}$ (and so these strata are $G$-orbits)
except for the open stratum $X_n\setminus X_{n-1}$ in the
case where $m=n$. In the latter case the $G$-orbits contained
in $X_n\setminus X_{n-1}$ are the level sets of the
determinant.

\bthm\label{transmat} Given two finite ordered collections $\cB$
and $\cB'$ of distinct matrices in $\Mat(n,m)$ of the same
cardinality, with the same sequence of ranks, and in the
case where $m=n$ with the same sequence of determinants, we can
simultaneously transform  $\cB$ into  $\cB'$ by means of an
element $g\in G$, where $G\subseteq\SAut(\Mat(n,m))$ is the
subgroup generated by all elementary one-parameter unipotent
subgroups. \ethm

Choosing in particular $\cB'$ consisting of diagonal matrices
we obtain a simultaneous diagonalization of the matrices in
$\cB$.

Theorem \ref{transmat} is an immediate consequence of Theorem
\ref{maincorstrat} and Lemma \ref{lemgeneric} below. To formulate
this lemma, we let $\cN$ be the saturated set of locally
nilpotent derivations on $X$ generated by all locally nilpotent
vector fields $\delta_{kl}$ and $\delta^{kl}$ ($k\ne l$) that is,
the set of all conjugates of these derivations along with their
replicas. The important observation is the following lemma.

\blem\label{lemgeneric}
$\cN$ has the orbit separation property on  $\Omega =X$.
\elem

\bproof In view of Lemma \ref{3.2} it suffices to show that the
derivations $\delta_{kl}$ and $\delta^{kl}$ have the orbit
separation property. Clearly it suffices to prove this for
$\delta^{kl}$. The action of the corresponding one-parameter
subgroup $H^{kl}=\exp(\kk \delta^{kl})$ on a matrix
$B=(b_1,\ldots, b_m)\in X$ with column vectors $b_1,\ldots, b_m\in
\kk^n$ is explicitely given by
$$
\exp(t\delta^{kl}).B= (b_1, \ldots, b_l+tb_k, \ldots, b_n),
$$
where $b_l+tb_k$ is the $l$th column
of the matrix on the right.
Thus the $H^{kl}$-orbit of $B$ has dimension one
if and only if $b_k\ne 0$. The functions
$$
B\mapsto b_{ij} \; (j\ne l)\and B\mapsto
\left|\begin{array}{cc}b_{ik} & b_{il} \\b_{jk} & b_{jl}
\end{array}\right|\; (i\ne j)
$$
on $\Mat(n,m)$ are $H^{kl}$-invariants that obviously
separate all $H^{kl}$-orbits of dimension one,
as the reader may easily verify.
\eproof

\subsection{The case of symmetric and skew-symmetric matrices}
\label{symasym}

We can apply the same reasoning to the varieties
$$
X=\Spec (\kk[T_{ij}]_{1\le i,j\le n}/(T_{ij}-T_{ji})_{1\le i,j\le n})
$$
of symmetric $n\times n$ matrices over $\kk$ and to the variety
$$
Y=\Spec (\kk[T_{ij}]_{1\le i,j\le n}/(T_{ij}+T_{ji})_{1\le i,j\le n})
$$
of skew symmetric matrices.
The group $\SL_n$ acts on both varieties via
$$
A.B=A BA^T\,,\mbox{where $A\in \SL_n$ and $B\in X (\in Y$, resp.)}.
$$
The subvariety $X_r$ of symmetric matrices of rank $r$
in $X$ is stabilized by this action, and also the determinant
of a matrix is preserved.
In the skew symmetric case again the subvarieties $Y_r$ of
matrices of rank $r$ are stabilized, and also the Pfaffian
$\Pf(B)$\footnote{We keep the usual convention that
the Pfaffian of a matrix of odd order equals zero. } of a matrix
$B\in Y$ is preserved.

By a well known theorem in linear
algebra the orbits of the $\SL_n$-action on $X$
are the subsets $X_r$ of matrices of rank $r$ in
$X$ for $r<n$, whereas for $r=n$ the orbits are
the level sets of the determinant. Similarly,
the orbits of the $\SL_n$-action on $Y$
are the subsets $Y_r$ for $r<n$, whereas for
$r=n$ the orbits are the level sets of the Pfaffian.

As in subsection 3.2 the elementary matrix $E_{ij}$ ($i\ne j$)
generates a one-parameter subgroup $H_{ij}=I_n+\kk E_{ij}$ that
acts on $X$ and on $Y$. The corresponding locally nilpotent vector
field will be denoted by $\delta_{ij}$.

Let $G_{\rm sym}\subseteq \SAut(X)$ and
$G_{\rm skew}\subseteq \SAut(Y)$, respectively,
be the subgroups generated by all $H_{ij}$
along with their replicas.
With these notations  the following results hold.

\bthm\label{transmatsym} Let $M_1, \ldots, M_k$ be a sequence of
pairwise distinct  symmetric  matrices of order $n\ge 2$ over
$\kk$. Assume that $M'_1, \ldots, M'_k$ is another such sequence
with
$$
\rk (M_i)=\rk(M_i')\ge 2 \and  \det (M_i)=\det (M_i')
\quad\forall i=1,\ldots k\,.
$$
Then there exists an automorphism $g\in G_{\rm sym}$ with
$g.M_i=M_i'$ for $i=1,\ldots k.$
\ethm

A similar result holds in the skew symmetric case.

\bthm\label{transmatskew} Let $M_1, \ldots, M_k$ be a sequence of
pairwise distinct skew-symmetric matrices of order $n\ge 2$ over
$\kk$. Assume that $M'_1, \ldots, M'_k$ is another such sequence
with
$$
\rk (M_i)=\rk(M_i') \and  \Pf (M_i)=\Pf (M_i')
\quad\forall i=1,\ldots k\,.
$$
Then there exists an automorphism $g\in G_{\rm skew}$ with
$g.M_i=M_i'$ for $i=1,\ldots k.$ \ethm

We give a sketch of the proof in the symmetric case only and leave
the skew-symmetric one to the reader. As in the case of generic
matrices (see Theorem \ref{transmat}) Theorem \ref{transmatsym} is
an immediate consequence of Theorem \ref{mainnewstrat} and Lemma
\ref{lemsym} below. In this lemma we let $\cN$ be the saturated
set of locally nilpotent derivations on $X$ generated by all
locally nilpotent vector fields $\delta_{kl}$.

\blem\label{lemsym}
$\cN$ has the orbit separation property on  $\Omega =X$.
\elem

\bproof In view of Lemma \ref{3.2} it suffices to show that the
derivations $\delta_{kl}$ have the orbit separation property. We
only treat the case $k<l$ the other one being similar. The action
of the corresponding one-parameter subgroup $H_{kl}=\exp(\kk
\delta_{kl})$ on a matrix $B\in X$ with entries $b_{ij}=b_{ji}$ is
explicitly given by $\exp(t\delta_{kl}).B= (b_{ij}')$, where
$$
b_{ij}'=b_{ij}\mbox{ if }i,j\ne k\,,\quad b'_{ki}=
b'_{ik}=b_{ik}+tb_{il} \mbox{ if } i\ne k\,,\and
b'_{kk}=b_{kk}+ 2t b_{kl}+t^2b_{ll}\,.
$$
Thus the $H_{kl}$-orbit of $B$ has dimension 0
if and only if $b_{il}= 0\,\forall i$. The functions
$$
B\mapsto b_{ij} \; (i,j\ne k)\,,\quad B\mapsto
\left|\begin{array}{cc}b_{ik} & b_{il} \\b_{jk} &
b_{jl}\end{array}\right|\; (i, j\ne k)\,,
\and
B\mapsto \left|\begin{array}{cc}b_{kk} & b_{kl}
\\b_{lk} & b_{ll}\end{array}\right|
$$
are  $H_{kl}$-invariants that are easily seen to separate all
$H_{kl}$-orbits of dimension one.
\eproof

\section{Tangential flexibility, interpolation by automorphisms, and
$\AA^1$-richness}\label{solupro}

\subsection{Flexibility of the tangent bundle}

We start with the following fact (see the Claim in the proof of
Corollary 2.8  in \cite{KK2}).

\blem
\label{5.10} Let $\p$ be a locally nilpotent vector field on
the affine $\kk$-scheme $X=\Spec A$ and let $p \in X$ be a point.
Assume that $f\in\ker \p$ is an invariant of $\p$ with $f(p)=0$.
If $\Phi=\exp(f\p)$ is the automorphism associated with the
locally nilpotent vector field $f \p$, then
\be\label{5.10a}
{d}_p\Phi (w) =w + {d}f (w) \p(p)\quad \mbox{ for all }w\in
T_pX\,. \ee
\elem

\bproof The tangent space $T_pX$ is the space of all derivations
$w:A\to \kk$ centered at $p$. For such a tangent vector $w$
its image ${d}\Phi(w)\in T_pX$ is the derivation \bals A\ni
g\mapsto w(\Phi(g))&= w\left(\sum_{i\ge 0}
\frac{f^i\p^i(g)}{i!}\right)
=\sum_{i\ge 0} \frac{1}{i!}w(f^i\p^i(g))\\
&=w(g)+w(f)\p(g)(p)\,, \eals as $f(p)=0$. Since by definition
$w(f)={d}f(w)$, the result follows.
\eproof

Now we can show the following result.

\bthm\label{5.15} Let $X$ be an irreducible affine algebraic
variety and let $G\subseteq \Aut(X)$ be a subgroup generated by a
saturated set $\cN$ of locally nilpotent vector fields. Assume
that $\cN$ satisfies the orbit separation property on a $G$-orbit
$O$. Then for each point $p\in O$, associating to an automorphism
$g\in G_{\cN,p}$ its tangent map  $dg(p)$ yields a representation
$$
\tau: G_{\cN,p}\lto \GL (T_pO)
\quad\mbox{with}\quad
\tau (G_{\cN,p})= \SL(T_pO)\,.
$$
\ethm

\bproof The assertion is trivially true if $\dim O=1$. Let us
assume for the rest of the proof that $\dim O\ge 2$. For any
one-parameter unipotent subgroup $H$ in $G_{\cN,p}$ the image
$\tau(H)$ is a subgroup of $\SL(T_pO)$. Hence also
$\tau(G_{\cN,p})\subseteq \SL(T_pO)$. Let us show the converse
inclusion.

According to Proposition \ref{1.5} there
are locally nilpotent vector fields
$\p_1, \ldots , \p_s\in\cN$
 spanning $T_xO$
at every point $x\in O$. Let $H_j=\exp(\kk \p_j)$ be the
one-parameter subgroup associated with $\p_j$. Using Remark
\ref{exanew}(3)  there are $H_j$-invariant open subsets
$U(H_j)\subseteq O$, $j=1,\ldots,s$ such that the geometric
quotients $\varrho_j:U(H_j)\to U(H_j)/H_j$
exist and satisfy the same properties as in \ref{exanew}(3).
In particular, the image $\vr_j(x)$ of a
generic point $x\in \bigcap_{j=1}^s U(H_j)$ is a smooth point
of  $U(H_j)/H_j$,
and $\vr_j$ has maximal rank at $x$.

We may assume that $\p_1(x),\ldots,\p_{m}(x)$, where $m=\dim
O\le s$, form a basis of $T_xO$. Hence for
$j,\mu\in\{1,\ldots, m\}$ with $\mu\ne j$ there exist
$\p_j$-invariant functions $f_{\mu j}$ on $X$ such that $f_{\mu
j}(x)=0$ and ${d}_x f_{\mu j}(\p_i(x)) =\delta_{\mu i}$. Consider
the automorphism $\Phi^t_{\mu j}=\exp(t\cdot f_{\mu j}\p_j)\in
G_{\cN,x}$ for $t \in \kk$. According to Lemma \ref{5.10} its
tangent map at $x$ is
$$
{d}_x\Phi^t_{\mu j}(\p_i(x))= \p_i(x) +t\cdot {d}_xf_{\mu j}
(\p_i(x))\cdot \p_j(x)= \p_i(x)  +t\delta_{\mu i}  \p_j(x)\,.
$$
Thus representing the elements in $\GL(T_xO)$ by matrices with
respect to the basis $\p_1(x),\ldots ,\p_m(x)$, the elements $
{d}_x\Phi^t_{\mu j}\in \GL(T_xO)$, $t\in \kk$, form just the
one-parameter unipotent subgroup generated by the elementary
matrix $E_{j\mu}$.  Since such one-parameter subgroups generate
$\SL(T_xO)$, the image of $G_{\cN,x}$ in $\GL(T_xO)$ contains
$\SL(T_xO)$ for a general point $x\in X$. Now the transitivity of
$G$ on $O$ implies that the same is true for every point $p\in O$.
\eproof

The following corollary is immediate.

\bcor\label{5.20} Under the assumptions of Theorem \ref{5.15}
for each point $p\in O$ we have
$$
\cN(p):=\{\p(p)\in T_pX: \p \in \cN\}=T_pO\,.
$$
In particular, the nilpotent cone $\LND_p(G)$ coincides with the
tangent space $T_p O$ for each $p\in O$. \ecor

\bproof Indeed, the group $G_{\cN,p}$ stabilizes $\cN(p)$ and for
$m=\dim O\ge 2$ the group $\SL(T_pO)$ acts transitively on
$T_pO\backslash\{0\}$. \eproof

\brem\label{nilcone} The last assertion in Corollary \ref{5.20}
does not hold any more for a general $\G_a$-generated subgroup
$G\subseteq\SAut(X)$ which is not generated by a saturated set of
locally nilpotent vector fields. For instance, if a semisimple
algebraic group $G$ acts on itself via left multiplications (i.e.,
$X=G$), then the cone $\LND_e(G)$ is just the usual nilpotent cone
in the Lie algebra $\Lie(G)=T_eX$, which is a proper subcone.
%To be more concrete, in the case $G=\SL_2$ the nilcone $\LND_e(G)$
%is just the quadratic cone in $\sl_2\cong\AA^3$ consisting of
%matrices with determinant $0$.
\erem

We also have the following  result on tangential
flexibility.\footnote{We are grateful to Adrien Dubouloz whose
observation allowed us to remove an inaccuracy in the original
formulation.}

\bcor\label{5.23} Let $\pi: E\to X$ be an irreducible and reduced
linear space\footnote{in the sense of \cite{EGA} Chap.\ II, 1.7.}
over a flexible variety $X$, which is over $X_\reg$ a vector
bundle. Assume that there is an action of $G:=\SAut(X)$ on $E$
such that the action of every 1-parameter subgroup of $G$ is
algebraic on $E$ and  $\pi$ is equivariant. Then the total space
$E$ is also flexible. In particular, the tangent bundle $TX$ and
all its tensor bundles $E=(TX)^{\otimes a}\otimes (T^*X)^{\otimes
b}$ are flexible. \ecor

\bproof It suffices to check that the special automorphism group
$G'=\SAut(E)$ acts transitively on $E_\reg=\pi^{-1}(X_\reg)$. By
our assumptions $G$ can be considered as an algebraically
generated subgroup of $G'$.  Since $X$ is flexible and $\pi$ is
equivariant, this subgroup acts transitively on the set of fibers
of $E_\reg\to X_\reg$. Moreover, $X$ being affine for any point
$e\in E_\reg$ there is a section $V:X\to E$ with $V(\pi(e))=e$.
This section generates a $\G_a$-action $w\mapsto w+tV(\pi(w))$.
Hence $G'$ acts transitively on every fiber of $E$ over a regular
point, and the result follows.\eproof

\bcor\label{5.22} Let $X$ be a flexible irreducible affine
variety of dimension $\ge 2$. Consider the special automorphism
group $G=\SAut(TX)$ of the tangent bundle $TX$, and let
$Z\subseteq TX$ be the zero section. Then the group $G_Z$ acts
infinitely transitively on $TX_{\reg}\backslash Z$. \ecor

\bproof The special automorphism group $\SAut(X)$ induces a
 $\G_a$-generated subgroup $\tilde G\subseteq G$ acting on $TX_{\reg}$.
Since $X$ is flexible this action is transitive on the zero
section, hence also on the set of fibers of $TX\to X$ over
$X_{\reg}$. On the other hand, by Theorem \ref{5.15} the
stationary subgroup $\tilde G_p$ of a given point $p \in X_{\reg}$
acts on $T_pX$ as $\SL(T_pX)$. Since $\dim T_pX>1$, it acts
transitively off the origin. Finally the action of $\tilde G$ on
$TX_{\reg}$ is transitive off the zero section. Hence by Theorem
\ref{maincor} the group $G_Z$, being $\G_a$-generated and
generated by a saturated set of locally nilpotent derivations,
acts infinitely transitively on $TX_{\reg}\backslash Z$. \eproof

\bsit\label{5.21} For later use let us mention the following
slightly more general version of Theorem \ref{5.15}. For a finite
subset $Z\subseteq X$ and $p\in O$  we let $\cN^M_{p,Z}\subseteq
\cN$ denote the set of all locally nilpotent vector fields $\p\in
\cN$ such that $\p$ has a zero at $p$ and a zero of order $\ge
M+1$ at all points of $Z\backslash \{p\}$. Let further $G_{p,Z}^M$
be the subgroup of $G$ generated by all exponentials of elements
in $\cN_{p,Z}^M$.  Replacing in Theorem \ref{5.15} $G_{\cN,p}$ by
$G_{p,Z}^M$ the following result holds. \esit

\bprop\label{5.15ext}
If $\dim O\ge 2$ then the image of the group $G_{p,Z}^M$ in $\GL(T_pO)$
coincides with $\SL(T_pO)$.
\eprop

\bproof With the notation as in the proof of {\em loc.cit.}, by
infinite transitivity (see Theorem \ref{mainnewstrat}) it suffices
to show the assertion for the case that $x=p$ is general and $Z$
consists of general points. Under this assumption we can find
$\p_j$-invariant functions $h_j$ with $h_j(x)=1$ which vanish in
all points of $Z\backslash\{p\}$. Replacing in the proof of
\ref{5.15} $f_{\mu j}$ by $h_j^{M+1}f_{\mu j}$, the automorphisms
$\Phi^t_{\mu j}$ are the identity up to order $M$ at the points of
$Z$ and remain unchanged at $x$. Now the same arguments as before
give the conclusion.\eproof

Let further $G^M_{Z}$ have the same meaning as $G^M_{p,Z}$ above,
but without any constraint imposed at $p$. That is, $G^M_{Z}$ is
the subgroup of $G$ generated by the saturated set $\cN_{Z}^M$ of
locally nilpotent vector fields vanishing to order $M+1\ge 1$ at
all points of $Z$. Then the same argument as before proves the
following proposition.

\bprop
\label{5.15ext2}
Every point $p\in O\backslash Z$ is $G^M_{Z}$-flexible,
hence $G^M_{Z}.p=O\backslash Z$.
\eprop

\subsection{Prescribed jets of automorphisms}

Let us start with the following standard fact (see Proposition 6.4.  in
\cite{KK1}, cf.\ also Theorem \ref{5.15}). Recall
that a {\em volume form} $\omega$ on a smooth algebraic variety
$X$ is a nowhere vanishing top-dimensional regular form on $X$; it
does exist if and only if $K_X=0$ in $\Pic(X)$.

\blem\label{5.30} If $X$ is an irreducible affine algebraic
variety and $\omega\in \Omega^n_X$ a volume form on $X_\reg$, then
$\omega$ is preserved under every automorphism $g\in \SAut(X)$.
\elem

\bproof It suffices to show that for every locally nilpotent
vector field $\partial$ the form $\omega$ is invariant under an
automorphism of $H=\exp(\kk \p)$. If $h_t=\exp(t\p)$ then for
every $x\in X_\reg$ the pullback $h_t^*(\omega)(x)$ is a multiple
of $\omega(x)$, i.e.\ $h_t^*(\omega)(x)=f(x,t)\omega(x)$, where
$f(x,t)\ne 0$ for all $x,t$. For a fixed $x$ the function $f(x,t)$
is thus  a polynomial in one variable without zero. Hence $f$ is
independent of $t$ and is equal to $f(x,0)=1$. \eproof

\bsit\label{5.31} We adopt the following notation and
assumptions. If $\varphi:X\to X$ is a morphism then its $m$-jet
$j^m_p\varphi$ at $p\in X$ can be regarded as a map of
$\kk$-algebras
$$
j^m_p\varphi:\cO_{X,\varphi(p)}/\fm_{\varphi(p)}^{m+1}
\lto \cO_{X,p}/\fm_p^{m+1}\,,
$$
where $\cO_{X,x}$ denotes the local ring at a point $x\in X$ and
$\fm_x$ its maximal ideal.

We assume in the sequel that $p\in X_{\reg}$ is a regular point
and $\varphi(p)=p$.  Letting $A_m=\cO_{X,p}/\fm_p^{m+1}$ the
$m$-jet of $\varphi$ yields a map of $k$-algebras
$$
j^m\varphi= j^m_p\varphi:A_m\lto A_m\,,
$$
which stabilizes the maximal ideal $\fm$ of $A_m$ and all of its
powers $\fm^k$.

For $m\ge 1$ we let $\Aut_{m-1}(A_m)$ denote the set of
$\kk$-algebra isomorphisms $f:A_m\to A_m$ with $f\equiv
\id\mod\fm^{m}$. For every $f\in \Aut_{m-1}(A_m)$ the map $f-\id$
sends $A_m$ into $\fm^m$ and vanishes on the constants $\kk$. As
it vanishes as well on $\fm^2$ it induces a $\kk$-linear map
$$
\psi_f:\fm/\fm^2\lto \fm^m =\fm^m_p/\fm_p^{m+1} \,.
$$
Note that  $\fm^m$ is naturally isomorphic to the $m$th symmetric power
$S^mV$ of the $\kk$-vector space $V=\fm/\fm^2$. For every $m\ge
1$ our construction yields a map \be\label{5.30a}
\psi:\Aut_{m-1}(A_m) \lto \Hom_\kk(\fm/\fm^2, \fm^m )\cong
V^\vee\otimes S^mV\,, \ee where $V^\vee$ stands for the dual
module of $V$. For $m=1$ this map associates to $f=j^1\varphi$
just the cotangent map $d\varphi(0)^\vee$.

In terms of local coordinates this construction can be interpreted
as follows. The $\kk$-algebra $A_m$ is isomorphic to the quotient
$A/\fm_A^{m+1}$, where $A=\kk[[x_1,\ldots, x_n]]$ is the
$\kk$-algebra of  formal power series and $\fm_A$ is its maximal
ideal.  Any map $f\in \Aut_{m-1}(A_m)$ is represented by an
$m$-jet of an $n$-tuple of power series $F=(F_1,\ldots, F_n)\in
A^n$ with $F_i\equiv x_i\mod\fm_A^{m}$. Clearly for any $m\ge 1$
the $m$-form $\psi_f$ corresponds to the $m$th order term of $F$.
\esit

With this notation we have the following lemma.

\blem\label{5.32}
\bnum[(a)]
\item For every $m\ge 1$ the map $\psi$ in
\eqref{5.30a} is bijective.
\item If $m=1$ then $\psi_{f\circ
g}=\psi_f\circ\psi_g$ while for $m\ge 2$ we have $\psi_{f\circ
g}=\psi_f+\psi_g$.
\item If $\p$ is a locally
nilpotent vector field on $X$ with a zero of order $m\ge 2$ at $p$
then $\psi_{\exp(t\p)}=t \psi_{\exp(\p)}$. \enum\elem

\bproof (a) is immediate using the coordinate description above.

(b) is easy and can be left to the reader. To deduce (c) we note
that $\exp(t\p)\in\SAut(X)$ induces the map $\id+t\hat\p\in
\Aut_{m-1}(A_m)$, where
$\hat \p$ denotes the derivation on $A_m$ induced by $\p$. Hence
$\psi_{\exp(t\p)}=t\hat\p$, proving (c). \eproof

An $n$-tuple $F=(F_1,\ldots, F_n)\in A^n$ as in \ref{5.31}
representing an $m$-jet $f=j^mF\in \Aut_{m-1}(A_m)$ preserves a
volume form $\omega$ on $X_{\reg}$ (or on $(X,p)$) if and only if
the Jacobian determinant $J_F$ of $F$ is equal to $1$. Modulo
$\fm^m$ this determinant depends only on $f$ and not on the
representative $F$ of $f$. Hence we can set $J_ f:= J_F\mod
\fm^m$.\footnote{However $J_f$ is not an element in $A_m$ since it
is  not well defined modulo $\fm^{m+1}$.} We say in the sequel
that an $m$-jet $f\in \Aut_{m-1}(A_m)$ with $m\ge 1$ preserves a
volume form  if $J_f\equiv 1\mod\fm^m$. The latter condition can
be detected in terms of $\psi_f$ as follows.

\blem\label{5.32a}
\bnum[(a)] \item
If $m=1$ then $f\in \Aut(A_1)$ preserves a
volume form if and only if  $\psi_f\in
\SL(V)$.

 In case $m\ge 2$ the map  $f\in \Aut_{m-1}(A_m)$
preserves a volume form if and only if $\psi_f$ is in the kernel
of the natural contraction map
\bals
\kappa_m: \Hom_\kk(V,S^mV)\cong V^\vee\otimes S^mV&\lto S^{m-1}V\,, \\
\lambda\otimes v_1\cdot \ldots \cdot v_m& \longmapsto
\sum_{\mu=1}^m \lambda(v_\mu)\cdot v_1\cdot\ldots\cdot \hat
v_\mu\cdot\ldots\cdot v_m.
\eals

\item $\ker\kappa_m$ is an irreducible $\SL_n(V)$-module for all
$m\ge 1$.
\enum \elem

\bproof In case $m=1$ (a) is immediate. Suppose that $m\ge 2$. If
$f= \id +f_m\mod \fm^{m+1}$ with an $n$-tuple of $m$-forms
$f_m=(f_{m1},\ldots, f_{mn})$, then  $J_f$ is easily seen to be
equal to
$$
1+\div f_m=1+\frac{\p f_{m 1}}{\p x_1}+
\ldots + \frac{\p f_{m n}}{\p x_n}
\quad\mod\fm^m\,,
$$
where $\div f_m$ is the divergence of $f_m$.
Thus $J_f\equiv 1\mod \fm^m$ if and only if $\div f_m=0$.
Writing $f_m\in V^\vee\otimes S^mV$ as
$f_m=\sum_{i=1}^n \frac{\p}{\p x_i}\otimes f_{m i}$
the element $\div f_m$ in $ S^{m-1}V$
corresponds just to the contraction $\kappa_m(f_m)$,
proving (a).

(b) is a standard fact in representation theory,
see e.g.\ \cite[\S IX.10.2]{Pr}.
\eproof

Now we can state our main result in this subsection.

\bthm\label{5.33} Let $X$ be an irreducible affine algebraic
variety of dimension $n\ge 2$ equipped with an algebraic  volume
form $\omega$ defined on $X_\reg$, and let $G\subseteq \SAut(X)$
be a subgroup  generated by a saturated set  $\cN$ of locally
nilpotent derivations. If $G$ acts  on $X$ with an open orbit $O$,
then for every $m\ge 0$ and every finite subset $Z\subseteq O$
there exists an automorphism $g \in G$ with prescribed $m$-jets
$j_p^m$ at the points $p\in Z$, provided these jets preserve
$\omega$ and inject $Z$ into $O$. \ethm

The proof will be reduced to the following lemma.

\blem\label{lemjet} With the notation and assumptions of Theorem
\ref{5.15}, suppose that $j^m_p$ is an $m$-jet of an automorphism
at a given point $p\in Z$, which is the identity up to order
$m-1\ge 0$. Then for every $M>0$ there is an automorphism $g\in G$
such that its $m$-jet at $p$ is $j_p^m$ while its $M$-jet at each
other point $q\ne p$ of $Z$ is the identity. \elem

Before proving Lemma \ref{lemjet} let us show how Theorem
\ref{5.33} follows.

\bproof[Proof of Theorem \ref{5.33}] We proceed by induction on
$m$. If $m=0$ the assertion follows from the fact that $G$ acts
infinitely transitively on $O$. For the induction step suppose
that we have an automorphism $g\in G$ with the prescribed jets up
to order $m-1\ge 0$. Thus the $m$-jets $j_p^{\prime m}=j_p^m \circ
g^{-1}$ are up to order $m-1$ the identity at every point $p\in
Z$. If we find an automorphism $h\in G$ with $m$-jet equal to
$j_p^{\prime m}$ for all $p\in Z$, then obviously the automorphism
$h\circ g$ has the desired properties.

Thus replacing $j_p^{m}$ by $j_p^{\prime m}$ we are reduced
to show the assertion in the case that for all $p\in Z$
the $m$-jets $j_p^{m}$ are the identity up to order $m-1$,
where $m\ge 1$.

Applying Lemma \ref{lemjet}, for every point $p\in Z$
there is an automorphism $g_p\in G$ whose $m$-jet at $p$
is the given one while its $m$-jets at all other points
$q\in Z\backslash \{p\}$ are the identity.
Obviously then the composition (in arbitrary order)
$g=\prod_{p\in Z}g_p$ will have the required properties.
\eproof

\bproof[Proof of Lemma \ref{lemjet}]
In the case $m=1$ the assertion follows
from Theorem \ref{5.15} and
Proposition \ref{5.15ext}. So we may assume for
the rest of the proof that $m\ge 2$.

Consider the set $\cN_{mp,Z}^M$ of all locally nilpotent
derivations in $\cN$ with a zero of order $m$ at $p$ and of order
$M+1$ at all other points $q\in Z\backslash\{p\}$. Let $G^M_{mp,Z}$ be
the subgroup of $G$ generated by the exponentials of elements in
$\cN_{mp,Z}^M$ so that an automorphism in $G^M_{mp,Z}$ is the
identity up to order $(m-1)$ at $p$ and up to order $M$ at all
other points $q\in Z\backslash\{p\}$. With the notation as
introduced in \ref{5.31} let us consider the composed map
$$
\Psi: G^M_{mp,Z} \lto \Aut_{m-1}(A_m)\stackrel{\psi}{\lto}
\Hom_\kk(V, S^mV)\,,
$$
where  $\psi$ is as in (\ref{5.30a})  and the first arrow assigns
to an automorphism its $m$-jet at $p$. Using Lemma \ref{5.32a}(a)
it suffices to show that $\Psi$ maps $G^M_{mp,Z}$
surjectively onto the subspace $\ker \kappa_m$.

The group $G^M_{mp,Z}$ is generated by exponentials of vector
fields in $\cN^M_{mp,Z}$. Thus using Lemma \ref{5.32}(b), (c) the
image $\im(\Psi)$ of $\Psi$ is a linear subspace of $\Hom_\kk(V,
S^mV)$. We claim that this subspace is nonzero.

Indeed, consider a vector field $\p\in \cN$ with $\p(p)\ne 0$ and
the one-parameter subgroup $H=\exp(\kk \p)$. According to
Remark \ref{exanew}(3) there is an open dense $H$-invariant subset
$U(H)\subseteq O$ which admits a quasi-affine geometric quotient
$U(H)/ H$ with the same properties as in \ref{exanew}(3). By
infinite transitivity of the action of $G$ on $O$ we may assume
that $Z\subseteq U(H)$ is such that the image of $Z$ in the
quotient $U(H)/ H$ is contained in the regular part of $U(H)/ H$,
has the same cardinality as $Z$, and the projection $U(H)\to U(H)/
H$ is smooth in the points of $Z$. Thus we can find a regular
$H$-invariant function $f$ on $X$ with a simple zero at $p$, and
another such function $h$ with $h(p)=1$ and $h(q)=0$ for all $q\in
Z\backslash\{p\}$. Replacing $f$ by $h^{M+1}f$ we may assume that
$f$ has a zero of order $\ge M+1$ at all points of
$Z\backslash\{p\}$ and a simple zero at $p$. Then $g=\exp(f^m\p)$
is an automorphism in $G^M_{mp,Z}$ with $\Psi(g)=f^m \hat\p\ne 0$,
where $\hat \p$ is the derivation of $A_m$ induced by $\p$ (cf.\
Lemma \ref{5.32}(c) and its proof). This proves the claim.

The group $G_{p,Z}^M$ acts on $G_{mp,Z}^M$ by conjugation
$g.h=g\circ h\circ g^{-1}$, where $g\in G_{p,Z}^M$ and $h\in
G_{mp,Z}^M$. If we write $h=\id +h_m\mod \fm^{m+1}$ with a map
$h_m\in \Hom_\kk(V, S^mV)$ then  $g.h=\id +g\circ h_m \circ
g^{-1}\mod\fm^{m+1}$. The map $g$ induces
the cotangent map $(d_pg)^\vee$
on $V=(T_pX)^\vee$ and its $m$th symmetric power
$S^m((d_pg)^\vee)$ on $S^mV$.
Hence there is a commutative diagram 
\bdi
G_{p,Z}^M\times G_{mp,Z}^M &\rTo & G_{mp,Z}^M\\
\dTo<{d_p^\vee-\times \Psi} &&\dTo>\Psi \\
\SL(V)\times \Hom_\kk(V,S^mV) &\rTo & \Hom_\kk(V,S^mV)\,, 
\edi
\noindent
where the lower horizontal map is induced by
the standard
representation of $\SL(V)$ on $S^mV$. Since the map $G_{p,Z}^M\to \SL(V)$ is
surjective (see Theorem \ref{5.15} and Remark \ref{5.21}), the
image $\im(\Psi)$ of $\Psi$ is a non-trivial $\SL(V)$-module. By Lemma
\ref{5.32} this representation is contained in the kernel of the
contraction map $\kappa_m$. Since the latter kernel is irreducible
(see Lemma \ref{5.32a}(b)), it follows that
$\im(\Psi)=\ker\kappa_m$, as required. \eproof

\brem\label{5.24} If in the situation of Theorem \ref{5.33} each
of the jets $j^m_p$, $p\in Z$, fixes the point $p$ and preserves a
volume form,\footnote{Note that this is a local condition, see the
discussion before Lemma \ref{5.32a}.} then the conclusion of
Theorem \ref{5.33} remains valid without the requirement that
there is a global volume form on $X_\reg$. \erem

\brem\label{5.24bis} If  $X_\reg$ does not admit a  global volume
form i.e., $K_{X_\reg}\neq 0$, one can still formulate a necessary
condition for interpolation of jets by an automorphism from a
$\G_a$-generated group $G$, namely in terms of the `volume form
monodromy' of $G$. To define it we fix   a volume form $\omega_x$
on the tangent space $T_xX$ at some point $x\in X_\reg$, and
consider the stabilizer subgroup $G_x\subseteq G$. Every element
$g\in G_x$ transforms $\omega_x$ into $\chi_x(g)\cdot \omega_x$,
where $\chi_x(g) \in \G_m=\G_m(\kk)$. The map
$$
\chi_x: G_x\lto\G_m
$$
is then a character on $G_x$ which equals $1$ on $G_{\cN,x}$, see
Theorem \ref{5.15}. If $y\in X$ is a second point and $h\in G$ is
an automorphism with $h.x=y$ then $hG_xh^{-1}=G_y$ and $h$
transforms $\omega_y$ into $\omega_x$. Hence
$\chi_y(hgh^{-1})=\chi_x(g)$ for all $g\in G_x$. In particular the
image of $\chi_x$ forms a subgroup $\Gamma$ of $\G_m$ independent
of $x\in O$, which is called the {\em volume form monodromy} of
$G$.

The volume form monodromy can be a nontrivial discrete group as in
the case of $X=\SL_2/N(\T)$ and $G=\SAut(X)$, where
$N(\T)\subseteq\SL_2$ is the normalizer of the maximal torus
$\T\subseteq\SL_2$. Note that in this case
$X\simeq\PP^2\backslash C$, where $C$ is a smooth conic in
$\PP^2$, see \cite{Po5}. Using technique from \cite{KK1} one can
show that here $\Gamma =\{\pm 1\}$. \erem

\subsection{ $\AA^1$-richness}
An irreducible affine variety $X$ is called {\em
$\AA^1$-rich} if for every closed  subset $Y$ of codimension $\ge
2$ and every finite subset $Z\subset X \setminus Y$ there is a
regular map $\AA^1\to X$ whose image contains $Z$ and omits $Y$
\cite[\S 2]{KZ2}.

The following  corollary is immediate from the Transversality
Theorem \ref{5.40}. In particular this result shows that a
flexible irreducible affine variety is $\AA^1$-rich. In the
special case where $X=\AA_\C^n$ the latter also follows from the
Gromov-Winkelmann theorem, see \cite[\S 2, Proposition 1]{Wi}.

\bcor\label{5.41} Let as before $X$ be an irreducible affine
variety and let $G\subseteq \SAut (X)$ be a subgroup generated by
a saturated set $\cN$ of locally nilpotent derivations, which acts
with an open orbit $O\subseteq X$. Then for any finite subset
$Z\subseteq O$ and any closed subset $Y\subseteq X$ of codimension
$\ge 2$ with $Z\cap Y=\emptyset$ there is an orbit $C\cong\AA^1$
of a $\G_a$-action on $X$ which does not meet $Y$ and passes
through each point of $Z$ having prescribed jets at these points.
\ecor

\bproof In the case $\dim X=1$ this is trivially true. So assume
that $\dim X\ge 2$. Let $C$ be an orbit of a $\G_a$-action on $O$.
Since $G$ acts infinitely transitively on $O$ we may assume that
$Z\subseteq C$. By Theorem \ref{5.33} and Remark \ref{5.24},
applying an appropriate automorphism $g'\in G$ we may suppose as
well that $C$ has prescribed $m$-jets at the points of $Z$.
Indeed, the $m$-jets of automorphisms stabilizing a given point
$p\in O$ and having at this point the jacobian determinant equal
to $1$ modulo $\fm^m$ act transitively on the set of all $m$-jets
of smooth curves at $p$.

By Proposition \ref{5.15ext2}, using the notation as in
\ref{5.21}, the $\G_a$-generated group $G^m_{Z}$ acts transitively
in $O\backslash Z$. Applying now the Transversality Theorem
\ref{5.40}(b) to $G^m_{Z}$, $C\cap (O\backslash Z)$, and $Y\cap
(O\backslash Z)$ we can find an element $g\in G^m_{Z}$ with
$g.C\cap Y=\emptyset$. Thus the $\G_a$-orbit $g.C$ contains $Z$,
has the prescribed jets at the points of $Z$, and does not meet
$Y$. \eproof

We can deduce also
the following fact.

\bprop\label{5.60} Let $G\subseteq \SAut (X)$ be a subgroup
generated by a saturated set $\cN$ of locally nilpotent
derivations, which acts with an open orbit $O\subseteq X$. Then
for any closed subset $Y\subseteq O$ of codimension $\ge 2$ and
for any $m\in\N$ the group $G^m_{\cN,Y}$ as in Remark \ref{Gire}
acts with an orbit open in $X$. \eprop

\bproof According to Proposition \ref{1.5} there are locally
nilpotent vector fields $\p_1,\ldots,\p_s$ generating $T_pX$ for
all $p\in O$. Let $H_\sigma\subseteq G$ be the one-parameter
subgroup associated to $\p_\sigma$. By Remark \ref{exanew}(3)
for suitable open dense $H_\sigma$-invariant subsets $U(H_\sigma)$
in $O$ there are geometric quotients $U(H_\sigma)/ H_\sigma$ as in
\ref{exanew}(3). Using the same reasoning as in the proof of
Theorem \ref{2.9} there is an $H_\sigma$-invariant function
$f_\sigma\in\cO(X)$ vanishing on $\ol{H_\sigma. Y}$ and equal to $1$ at a
given general point $p\in U(H_\sigma)\backslash \ol{H_\sigma. Y}$.
Consequently the replica $\exp(\kk f^m_\sigma\p_\sigma)$
fixes the $m$th infinitesimal neighborhood of $Y$ in $X$ and moves
$p$ in direction $\p_\sigma (p)$. In other words, $p$ is a
$G^m_{\cN,Y}$-flexible point. Applying Corollary
\ref{gflex1}(a) the result follows. \eproof

\bprob\label{codimtwo}   Is it true that in the situation of
Proposition \ref{5.60} the group $G^m_{\cN, Y}$ acts
transitively on $O\backslash Y$?
 \eprob

In the case that $X=\AA^n_\C$ and $G= \SAut (X)$ the answer is affirmative,
see \cite{Wi}, \S 2, Proposition 1 and its proof.

\brems 1. \label{jal}  Every algebraic variety $X$ contains a divisor
$Y$ such that the logarithmic Kodaira dimension
$\bar\kappa(X\backslash Y)$ is $\ge 0$. In this case $X\backslash
Y$ cannot carry a $\G_a$-action and so $G_Y=\{\rm id\}$
although $X$ might be flexible. The simplest example of such a
situation is given by the hypersurface $Y=\{X_1\cdot\ldots\cdot
X_n=0\}$ in $X=\AA^n$, see also \cite{Jel}.

2. \label{newre10} If the group $\SAut(X)$ of an irreducible
normal affine surface $X$ acts with an orbit $O$ open in $X$ then
$X$ is a Gizatullin surface  and $Y=X\setminus O$ is a finite set,
see \cite[II, Theorem 3]{Gi}, \cite{Du}. Such a surface $X$ is
usually non-$\Q$-factorial. In higher dimensions this complement
may contain a divisor (see Example \ref{singloc} below). However,
such examples cannot exist if $X$ is $\Q$-factorial (see Corollary
\ref{5.80} below). \erems

We need the following auxiliary result.

\bprop\label{5.70} Let $X$ be an irreducible normal affine
variety, and let $G$ be a $\G_a$-generated subgroup of $\SAut(X)$
acting on $X$ with an open orbit $O\subseteq X$. If the complement
$X\backslash O$ contains a divisor $D$ then $D$ generates a
nonzero element $[D]$ in the divisor class group $\Cl(X)_\Q =\Cl
(X) \otimes_\Z \Q$ over $\Q$. \eprop

\bproof
Assume to the contrary that $[D]=0$ in $\Cl(X)_\Q$.
Then there is a function $f$ on $X$ with $D=\V(f)$ set
theoretically. For every one dimensional unipotent subgroup $H\subseteq G$
and $x\in O$ the function $f|H.x$ is a polynomial on
$H.x\cong\kk$. As $H.x\subseteq O$ and so $D\cap H=\emptyset$,
this polynomial has no zero and so is constant equal to $a:=f(x)$. Hence
$H.x$ is contained in the level set  $f^{-1}(a)$ of $f$. Since $G$ is
generated by such subgroups, the whole orbit $O=G.x$ is
contained in $f^{-1}(a)$ and so it cannot be open, a
contradiction.
\eproof

\bcor\label{5.80} Let $G$ be a $\G_a$-generated subgroup of
$\SAut(X)$. If $X$ is $\Q$-factorial and a closed subset
$Y\subseteq X$ contains a divisor, then the group $G_{\cN,Y}$ has
no open orbit. \ecor

\section{Some applications}\label{appli}

\subsection{Unirationality, flexibility,
and triviality of the Makar-Limanov invariant}\label{fv.12.26.11}

Recall \cite{Fre} that the {\em Makar-Limanov invariant} $\ML(X)$
of an affine variety $X$ is the intersection of the kernels of all
locally nilpotent derivations on $X$. In other words $\ML(X)$ is
the subalgebra  of the algebra $\cO(X)$ consisting of all
$\SAut(X)$-invariants. Similarly \cite{Lie1} the {\em field
Makar-Limanov invariant} $\FML(X)$ is defined as the subfield of
$\cM er(X)$ which consists of all rational
$\SAut(X)$-invariants. If it is trivial i.e., if $\FML(X)=\kk$
then so is $\ML(X)$, while the converse is not true in general,
see Example \ref{rm1}(2) below. The next proposition confirms, in
particular, Conjecture 5.3 in \cite{Lie1} (cf.\ also
\cite{Po2}).

\bprop\label{prep} An irreducible affine variety $X$
possesses a flexible point if and only if the group $\SAut(X)$
acts on $X$ with an open orbit, if and only if the field
Makar-Limanov invariant $\FML(X)$ is trivial. In the latter case
$X$ is unirational. \eprop

\bproof The first equivalence follows from Corollary
\ref{gflex1}(a) and the second from Corollary \ref{RTcor}. As for
the last assertion, see the next remark.
%Remark \ref{a1con} below.
\eproof

\brem\label{a1con} As follows from Proposition \ref{1.1}(b) for
every $G$-orbit $O$ of a $\G_a$-generated group
$G\subseteq\SAut(X)$ there is a surjective morphism $\AA^s \to O$.
Hence any two points in $O$ are contained in the image of a
morphism $\AA^1\to O$. In particular $O$ is {\em
$\AA^1$-connected} in the sense of \cite[\S 6.2]{KK2}. \erem

\bexas\label{rm1} 1. Flexibility implies  neither rationality nor
stable rationality. Indeed, there exists a finite subgroup
$F\subseteq\SL(n,\C)$, where $n\ge 4$, such that the smooth
unirational affine variety $X=\SL(n,\C)/F$ is not stably rational,
see \cite[Example 1.22]{Po2}. However, by Proposition \ref{afhom}
below $X$ is flexible and the group $\SAut(X)$ acts infinitely
transitively on $X$.

2.  There are non-unirational affine threefolds $X$ with
$\ML(X)=\kk$ birationally equivalent to $C\times \AA^2$, where $C$
is a curve of genus $g\ge 1$, see \cite[\S 4.2]{Lie}. For such a
threefold $X$ the general $\SAut (X)$-orbits have dimension
two, the field Makar-Limanov invariant $\FML(X)$ is non-trivial,
and there is no flexible point in $X$.
\eexas

\subsection{Flexible quasihomogeneous varieties}
An important class of flexible algebraic varieties consists of
homogeneous spaces of semisimple algebraic groups. More
generally, the following hold (cf.\ \cite[\S 1.1]{Po2}).

\bprop\label{afhom} Let  $G$ be a connected affine algebraic group
without non-trivial characters,
and let $H$ be a closed subgroup of
$G$. Then the homogeneous space $G/H$ is flexible. In particular,
if $G/H$ is affine of dimension $n\ge 2$ then the group $\SAut
(G/H)$ acts infinitely transitively on $G/H$.\eprop

\begin{proof}
The image of $G$  in $\SAut(G/H)$ is a $\G_a$-generated
subgroup (see Example \ref{fv.12.26.11a} (2)). Thus the group
$\SAut(G/H)$ acts on the quotient $G/H$ transitively and $G/H$ is
flexible; see Proposition 1.1 in \cite{AKZ}. The second assertion
follows from the first one in view of Theorem \ref{mthm} and
Corollary \ref{tr-fl}.
\end{proof}

The following problem arises.

\noindent \bprob\label{prohom} {\em Characterize flexible
varieties among affine varieties admitting an action of a
semisimple algebraic group with a dense open orbit.}\eprob

For instance, if such a quasihomogeneous variety is smooth then in
fact it is flexible. In the particular case $G=\SL_2$ this was
actually established  in \cite[III]{Po5}, where we  borrowed the
idea of the proof of the following theorem.

\bthm\label{quasihom} Suppose that a  connected semisimple
algebraic group $G$ acts on a smooth irreducible affine
variety $X=\Spec A$ with an open orbit. Then $X$ is homogeneous
with respect to a connected affine algebraic group $\tilde
G\supseteq G$ without non-trivial characters. In particular, $X$
is flexible. \ethm

\bproof Since by our assumption $A^G=\kk$, the variety $X$
contains a unique closed $G$-orbit $Z\subseteq X$ and the
stabilizer of a point on this orbit is a reductive subgroup $H$ of
the group $G$ (see Theorems 4.17 and 6.7 in \cite{PV}).
Moreover, it follows from Luna's \'Etale Slice Theorem that there
is a finite dimensional rational $H$-module $W$ such that the
variety $X$ is $G$-equivariantly isomorphic to the total space
$G\times^H W$ of the homogeneous vector bundle over $G/H$ with the
fiber $W$, see Theorem~6.7 in \cite{PV}.

%Moreover, it follows from Luna's \'Etale Slice Theorem that the
%variety $X$ is $G$-equivariantly isomorphic to the total space of
%the homogeneous fiber bundle $(G\times  W)/H \to G/H$, where $W$
%is a rational $H$-module and where $H$ acts on $G\times W$ via
%\be\label{acts} h.(g,w)=(gh^{-1},h.w)\,, \ee see Theorem~6.7 in
%\cite{PV}.

According to \cite{BBHM} there exists a finite dimensional
$G$-module $V$ such that $V=W\oplus W'$, where $W'\subseteq V$ is
a complementary $H$-submodule. Letting
$$\tilde G=G\ltimes V\quad\mbox{with}\quad (g_1,v_1).(g_2,v_2)=
(g_1g_2,g_2^{-1}v_1+v_2)\,,$$ $\tilde H=H\ltimes W'$,  and $\tilde
H_0=\{e\}\ltimes W'$
%, where $H$ acts on $G\times W$ via (\ref{acts}),
we can identify ${\tilde G}/\tilde H_0$ and $G\times W$ as
$H$-varieties. Since the subgroup $H\subseteq G$ normalizes
$\tilde H_0$ in $\tilde G$ it acts  $\tilde G$-equivariantly on
the right on ${\tilde G}/\tilde H_0$. The latter fact can be used
to deduce the isomorphisms of varieties
$${\tilde G}/\tilde H\cong ({\tilde G}/\tilde H_0)/H
%\cong (G\times W)/H
\cong X\,.
$$
By Proposition \ref{afhom} $X$ is flexible being a homogeneous
variety  of a connected affine algebraic group ${\tilde G}$
without non-trivial characters (indeed, $\tilde
G=(G,0)\cdot(e,V)$, where both groups admit no non-trivial
characters). Now the proof is completed. \eproof

In the next theorem we provide a complete solution of Problem
\ref{prohom} for $G=\SL_2:=\SL_2(\kk)$ and $X$ normal.

\bthm\label{sl2} Every normal irreducible affine variety $E$
admitting an $\SL_2$-action with an open orbit is flexible. \ethm

For a homogeneous affine variety $E=\SL_2/H$ the result follows
from Proposition  \ref{afhom}. The proof  in the general case
given below is based on a description of normal $\SL_2$-varieties
due to Popov \cite[I]{Po5} (see also \cite[Chapter III, \S 4]{Kr})
and a Cox ring $\SL_2$-construction due to Batyrev and Haddad
\cite{BH}. Recall \cite[I]{Po5} that every non-homogeneous 
normal affine $\SL_2$-threefold with an open orbit is uniquely
determined by a pair $(h,m)$, where $m$ is the order of the
generic isotropy group\footnote{Which is a cyclic group.} and $h=
p/q\in ]0, 1]\cap\Q$ is the so called {\em  height} of $X$.  The
$\SL_2$-threefold  with invariant $(h,m)$ is denoted by $E_{h,m}$.
Notice that $E_{h,m}$ is smooth for $h=1$ and singular for $h<1$.

Assuming in the sequel that $p$ and $q$ are coprime positive  integers
 we
let
\be\label{ab} a = m/k
\quad\mbox{and}\quad b = (q-p)/k ,\quad\mbox{where}\quad k
= \gcd(q-p,m)
\,.\ee Let $\mu_a = \langle \xi_a\rangle$
denote the cyclic group generated by
a primitive root of unity $\xi_a\in \G_m=\G_m(\kk)$ of degree $a$.
The $\SL_2$-variety $E_{h,m}$ is isomorphic to the categorical
quotient of the hypersurface $D_b \subseteq \AA^5$ with equation
\be\label{hyper} Y^b = X_1X_4 - X_2X_3\ee
modulo the diagonal action of the group $\G_m
\times\mu_a$ on $\AA^5=\Spec\kk[X_1,X_2,X_3,X_4,Y]$ via
$${\rm diag} (t^{-p}, t^{-p}, t^q, t^q,t^k)\times
{\rm diag} (\xi^{-1}, \xi^{-1}, \xi, \xi,1)\,, \quad t\in
\G_m,\,\xi\in\mu_a\,.$$ Here the $\SL_2$-action on $D_b$ is
induced by the trivial action on the coordinate $Y$, while
$\langle X_1,X_2\rangle$ and $\langle X_3,X_4\rangle$ are simple
$\SL_2$-modules. This $\SL_2$-action on $D_b$ commutes with the
$(\G_m \times\mu_a)$-action and so descends to the quotient. This
gives a simple and uniform description of all non-homogeneous
singular normal affine $\SL_2$-threefolds with an open orbit
$E_{h,m}$ via the Cox realization as the quotient of the spectrum
of the corresponding Cox ring by the action of the Neron-Severi
quasitorus, see \cite{BH}.

\bproof[Proof of Theorem \ref{sl2}.] Let $E$ be a non-homogeneous
singular  normal irreducible affine $\SL_2$-variety with an
open orbit. If $\dim E=2$ then $E$ is a toric surface, in fact a
Veronese cone, and the group $\SL_2$ is transitive off the vertex
(see \cite[II]{Po5} or, alternatively, Theorem 0.2 in \cite{AKZ}).
Now the assertion follows by Theorem \ref{mthm}.

If further $E$ as above is an $\SL_2$-threefold then according to
Popov's classification $E=E_{h,m}$ for some pair $(h,m)$.

In the case where $E=E_{h,m}$ is smooth that is, $h=1$ the result
follows from Theorem \ref{quasihom}.

In the case where $E=E_{h,m}$ is singular
i.e., $h=p/q<1$, there is a unique singular point, say, $Q\in E$.
The complement $E\backslash \{Q\}$ consists of two $\SL_2$-orbits
$O_1$ and $O_2$, where $O_1\cong \SL_2/\mu_m$ while $O_2\cong
\SL_2/U_{a(p+q)}$
has the isotropy subgroup $$U_{a(p+q)}=\left\{\left(%
\begin{array}{cc}
  \xi & \eta \\
  0 & \xi^{-1} \\
\end{array}%
\right) \,\vert\,\eta\in\kk,\,\xi^{a(p+q)}=1\right\}\,.$$
Consider the hypersurface $D_b \subseteq \AA^5$
as in (\ref{hyper}). We can realize $\AA^5$ as a matrix space:
$$\AA^5=\left\{(X,Y)\,\vert\, X=\left(%
\begin{array}{cc}
  X_1 & X_3 \\
  X_2 & X_4 \\
\end{array}%
\right), X_i, Y\in\AA^1\right\}\,.$$
Then according to \cite{BH} the 3-fold $E=E_{h,m}$
admits a realization
as the categorical quotient of $D_b$
 by the action of the group $G_m\times \mu_a$  via
 $$(t,\xi).(X,Y)=\left(\left(%
\begin{array}{cc}
  \xi^{-1}t^{-p}X_1 & \xi t^{q}X_3 \\
  \xi^{-1}t^{-p}X_2 & \xi t^{q}X_4 \\
\end{array}%
\right)\,, t^kY\right)\,.$$
This action commutes with the natural $\SL_2$-action on
$D_b$ given by
 $$A.(X,Y)=(AX,Y)\,.$$
Hence the  $\SL_2$-action  on $D_b$ descends to the quotient
$E=E_{h,m}$. The hypersurface $Z=\{Y=0\}$ in  $D_b$ is the inverse
image of the unique two dimensional $\SL_2$-orbit closure in
$E_{h,m}$. To show the transitivity (or the flexibility) of the
group $\SAut(X)$ in $E_{\reg}$ it suffices to find a locally
nilpotent derivation $\p$ of the algebra $\cO(D_b)$ with
$\p(Y)\neq 0$ which preserves the $(\Z\times\Z_a)$-bigrading on
$\cO(D_b)$ defined via
$$\deg X_1=\deg X_2=(-p,-\bar 1),\,\,\,\deg X_3=\deg X_4=(q,\bar 1),
\quad\mbox{and}\quad\deg Y=(k,\bar 0)\,.$$ Indeed, such a
derivation induces a locally nilpotent derivation on $\cO(E)$.
Since $\p(Y)\neq 0$ the restriction of the corresponding vector
field to the image $Z'$ of $Z$ in $E$ is nonzero and so the
points of $Z'$ with $\p\neq 0$ are flexible. By transitivity,
every point of $Z'\backslash\{Q\}$ is.

The variety $D_b$ can be regarded as a suspension
over $\AA^3=\Spec\kk[X_2,X_3,Y]$, see (\ref{SuS}). Namely,
$$D_b=\{X_1X_4=f(X_2,X_3,Y)\}\quad\mbox{where}\quad
f=X_2X_3+Y^b\,.$$ According to  \cite{AKZ} (see also
Lemma 3.3 in \cite[\S 5]{KZ1}) a desired bihomogeneous locally nilpotent
derivation $\p$ can be produced starting
with a locally nilpotent derivation $\delta\in\Der\kk[X_2,X_3,Y]$.
For instance, let $\delta$ be given by
$$\delta(X_2)=\delta(X_3)=0,\quad \delta(Y)=X_2^cX_3^d\,.$$ Then
$\p$ can be defined via
\be\label{der}\qquad \p(X_1)=\p(X_2)=\p(X_3)=0,\,\,\p(X_4)
=\delta(f)=bX_2^cX_3^dY^{b-1},\,\,\p(Y)=X_1X_2^cX_3^d\,\ee with $a,b$
as in (\ref{ab}) and with appropriate values of
the natural parameters $c,d$.
Such a derivation $\p$ preserves the
$(\Z\times\Z_a)$-bigrading\footnote{I.e. $\deg\p(Y)=\deg Y$
and $\deg\p(X_4)=\deg X_4$.} if and only if
$$
\begin{aligned} -p-cp+dq&=k\\ k(b-1)-cp+dq&=q\\
 -1-c +d&\equiv 0\mod a\\-c+d &\equiv 1\mod a\,.
\end{aligned}
$$ By virtue of (\ref{ab}) the second relation follows
from the first one, while the last
two are equivalent.
Letting $c=s-1$ we can rewrite the remaining relations as
\be\label{co}
\begin{aligned} dq-sp=k\\ s\equiv d\mod a\,.\end{aligned}\ee
Since $\gcd(p,q)=1$
the first equation
admits a solution $(d_0,s_0)$ in natural numbers.
For every $r\in \N$, the pair $(d_0+rp,s_0+rq)$
also represents such a solution.
The second relation in (\ref{co}) becomes
\be\label{cong0}
r(q-p)\equiv d_0-s_0\mod a\,.\ee
By (\ref{ab}) $k=\gcd(m,q-p)$,
hence $\gcd(k,p)=1$. The first equation  in (\ref{co}) written as
$$d_0(q-p)-p(s_0-d_0)=k$$ implies that $k\,\vert (s_0-d_0)$.

Let $l=\gcd(a,q-p)=\gcd(a,bk)$. Since $\gcd(a,b)=1$ then $l\,\vert k$
and so (\ref{cong0}) is equivalent to the congruence
$$r\cdot\frac{q-p}{l}\equiv \frac{d_0-s_0}{l}\mod \frac{a}{l}\,.$$
Since $\frac{q-p}{l}$ and $\frac{a}{l}$ are coprime
the latter congruence admits a solution, say, $r_0$.
Letting finally $$c=s_0+r_0q-1, \,d=d_0+r_0p$$ the
locally nilpotent derivation $\p$ as in (\ref{der})
becomes homogeneous of bidegree $(0,\bar 0)$, as needed.
Now the proof is completed.
\eproof

The question arises whether the smooth loci of
singular affine $\SL_2$-threefolds
are homogeneous as well, cf.\ Theorem \ref{quasihom}.
The answer is negative; the following
proposition gives a more precise information.

\bprop\label{sisl2} Let $E=E_{h,m}$, where $h=p/q<1$ with
$\gcd(p,q)=1$. The following conditions are equivalent: \bnum
\item[(i)] The $\SL_2$-action on $E$ extends to an action
of a bigger affine algebraic group $G$ on $E$ which is transitive
in $E_{\reg}$;\item[(ii)] The variety $E$ is toric;
\item[(iii)]
$(q-p)\,\vert m$ or, equivalently, $b=1$ in (\ref{ab}).
\enum \eprop

\bproof Implication (i)$\Rightarrow$(ii) follows from Theorem 1 in
\cite[III]{Po5}. According to this theorem, a normal affine
threefold $X$ with a unique singular point $Q$ which admits an
action of an affine algebraic group transitive on
$X\backslash\{Q\}$, is toric.

The equivalence (ii)$\Leftrightarrow$(iii) follows from the
results of \cite{BH} and \cite{Ga}. Let us show the remaining
implication (iii)$\Rightarrow$(i). If $b=1$ in  (\ref{ab}) then
$D_b\cong \AA^4=\Spec\kk[X_1,\ldots,X_4]$. Hence the toric variety
$E_{h,m}$ can be obtained as the quotient
$\AA^4/(G_m\times\mu_a)$, where the group $G_m\times\mu_a$ with
$a=m/(q-p)$ as in  (\ref{ab}) acts diagonally on $\AA^4$ via
\be\label{gmact}\qquad\,\,\,\, (X_1,X_2,X_3,X_4)\longmapsto
(\xi^{-1}t^{-p}X_1,\xi^{-1}t^{-p}X_2, \xi t^{q}X_3,\xi
t^{q}X_4),\quad (t,\xi)\in G_m\times\mu_a\,.\ee Consider the
action of the group $\SL_2\times \SL_2$ on  $\AA^4$
via $$(A_1,A_2).(X_1,X_2,X_3,X_4)=\left(A_1\left(%
\begin{array}{c}
  X_1 \\
  X_2 \\
\end{array}%
\right),\,\,A_2\left(%
\begin{array}{c}
  X_3 \\
  X_4 \\
\end{array}%
\right)\right)\,.$$ This action commutes with the
$(G_m\times\mu_a)$-action (\ref{gmact}) and so descends to the
quotient $E_{h,m}$. The induced $(\SL_2\times \SL_2)$-action on
the quotient $E_{h,m}$ is transitive in the complement of the
unique singular point $Q$. This yields (i). Now the proof is completed.
\eproof

\bcor\label{trcor} None of the non-toric affine threefolds
$E=E_{h,m}$ with $h<1$ admits an algebraic group action transitive in
$E_{\reg}$. However, the group $\SAut (E)$ acts infinitely transitively
in $E_{\reg}$.\ecor

Let us finish this subsection with an example of a  flexible
non-normal irreducible affine variety with singular locus of
codimension one.

\bexa\label{singloc} Consider the standard irreducible
representation of the group $\SL_2$ on the space of binary forms
of degree three
$$V=\left\langle X^3,X^2Y,XY^2,Y^3\right\rangle\,.$$
Restriction to the subvariety
$$E=\SL_2 . X^2Y\cup \SL_2 . X^3\cup \{0\}\subseteq V$$
of forms with zero discriminant yields a non-normal
$\SL_2$-embedding, see \cite{Kr}. Since for a hypersurface in a
smooth variety normality is equivalent to smoothness in
codimension one, the divisor $D=\SL_2 . X^3\cup \{0\}\subseteq E$
coincides with the singular locus $E_{\rm sing}$. The complement
$E_{\reg}=\SL_2 . X^2Y$ is the open $\SL_2$-orbit consisting of
all flexible points of $E$. Hence $E$ is flexible.

The normalization of $E$ is isomorphic to
$E_{\frac{1}{2},1}$. Indeed $m=1$ because the stabilizer in
$\SL_2$ of a general point in $E$ is trivial. On the other hand,
the order of the stabilizer of the two dimensional orbit equals
$p+q=3$, hence $p=1$ and $q=2$.

For any admissible pair $(h,m)$ the affine threefold $E_{h,m}$ is
a union of an open $\SL_2$-orbit $O$ and an invariant prime
divisor $Y$ \cite[I, Lemma~4 and Corollary~1]{Po5}. Choosing a
generating set of one-parameter unipotent subgroups of $\SL_2$ we
let $\cN$ be the corresponding saturated set of locally nilpotent
vector fields on $E$. Consider further the subgroup
$G\subseteq\SAut (E)$ generated by $\cN$. Clearly, $G$ again
stabilizes the divisor $Y$ and acts transitively on its complement
$O$. According to Theorem \ref{2.9} the group $G_Y$ also acts on
$E$ with an orbit $O$ whose complement $Y$ is a divisor. This
shows that the assumption of $\Q$-factoriality in Corollary
\ref{5.80} is essential.  \eexa

\section{Appendix: Holomorphic flexibility}\label{holomo}
In this appendix we extend the notion of a flexible affine variety
to the complex analytic setting (cf.\ \cite{For1}). We survey
relations between holomorphic flexibility, Gromov's spray and the
Andersen-Lempert theory. In particular, we show that every
flexible variety admits a Gromov spray. This provides a new wide
class of examples to which the Oka-Grauert-Gromov Principle can be
applied. We refer the reader to \cite{For3} and the survey
articles \cite[\S 3]{For2} and \cite{KK} for a more thorough
treatment and historical references.

\subsection{Oka-Grauert-Gromov Principle for flexible varieties}
\label{OGG} The following notions were introduced in \cite[\S
1.1.B]{Gro}.

\bdefi\label{spray} (i) Let $X$ be a complex manifold. A  {\em
dominating spray} on  $X$ is a holomorphic vector bundle $\rho : E
\to X$ together with a holomorphic map $s : E \to X$, such that
$s$ restricts to the identity on the zero section $Z$ while for
each $x \in Z\cong X$ the tangent map $d_xs$ sends the fiber
$E_x=\rho^{-1}(x)$ (viewed as a linear subspace of $T_x E$)
surjectively onto $T_x X$.

(ii) Let $h: X \to B$ be a surjective submersion
 of complex manifolds. We say that it admits
a {\em fiber dominating spray} if there is a holomorphic vector
bundle $E$ on $X$ together with a holomorphic map  $s : E \to X$
such that the restriction of $s$ to each fiber $h^{-1}(b),
\, b \in B$, yields a spray on this fiber. \edefi

In these terms, the Oka-Grauert-Gromov Principle can be stated as
follows.

\bthm\label{6.20} (\cite[\S 4.5]{Gro}) Let $h : X\to B$ be a
surjective submersion of Stein manifolds. If it admits
a fiber dominating spray then the following hold.
\bnum[(a)]
\item Any continuous section of $h$ is homotopic
to a holomorphic one; and
\item  any two holomorphic sections of $h$
that are homotopic via continuous sections are
also homotopic via holomorphic ones.
\enum
\ethm

Due to the following proposition, smooth affine algebraic
$G$-fibrations with flexible fibers are appropriate for applying
this principle (cf.\ \cite[3.4]{For2}, \cite{Gro}).

\bprop\label{6.10} \bnum[(a)]
\item Every flexible smooth irreducible affine algebraic variety
$X$ over $\C$ admits a dominating spray.
\item
Let $h : X \to B$ be a surjective submersion of smooth
irreducible affine algebraic varieties  over $\C$ such that for
some algebraically generated subgroup $G\subseteq\Aut (X)$ the
orbits of $G$ coincide with the fibers of $h$ \footnote{We say in
this case that $X$ is {\em $G$-flexible over} $B$.}. Then $X\to B$
admits a fiber dominating spray. \enum \eprop

\bproof
It suffices to show (b). Indeed, due to Corollary \ref{tr-fl},
(a) is a particular case of (b).

By Proposition \ref{1.5} there is a sequence of algebraic
subgroups $\cH=(H_1,\ldots,H_s)$ of $G$ such that the tangent
space to the orbit $G.x$ at each point $x \in X$ is spanned by the
tangent spaces at $x$ to the orbits $H_i .x$, $i=1,\ldots,s$. Let
$\exp: T_1(H_i)\to H_i$ be the exponential map. Letting $E=
X\times \prod_{i=1}^s T_1(H_i)$ be the trivial vector bundle over
$X$ we consider the morphism
$$
s: E \to X,\qquad (x, (h_1,\ldots,h_s) ) \longmapsto
\Phi_{\cH,x} (\exp h_1,\ldots,\exp h_s)\,,
$$
where $\Phi_{\cH,x}$ has the same meaning as in (\ref{1.1.a}).
This yields the desired dominating spray.
\eproof

To extend Proposition \ref{6.20} to the analytic setting we
introduce below the notions of holomorphic flexibility. Recall
that a holomorphic vector field on a complex manifold $X$ is {\em
completely integrable} if its phase flow defines a holomorphic
action  on $X$ of the additive group $\C_+=\G_a(\C)$.

\bdefis\label{6.30}  (i) We say that a Stein
space $X$ is {\em holomorphically flexible} if the completely
integrable holomorphic vector fields on $X$ span the tangent space
$T_xX$ at every smooth point of $X$.

(ii) Given a
 holomorphic submersion
$h: X \to B$ of Stein manifolds, we say that $X$
is {\em holomorphically flexible over} $B$ if the completely
integrable relative holomorphic vector fields on $X$ span the
relative tangent bundle of $X\to B$ at any point of $X$. In the
latter case each fiber $h^{-1}(b), \, b \in B$, is a
holomorphically flexible Stein manifold. \edefis

\brems\label{6.31} 1. The vector field $\delta=z\frac{d}{dz}$ on
$X=\C^*=\C\backslash\{0\}$ is completely integrable. However, the
derivation $\delta\in\Der(\cO(X))$ is not locally nilpotent. Hence
$X=\C^*$ is not flexible in the sense used in this paper, while it
is holomorphically flexible.

2. In the terminology of \cite{Va},
a complex manifold $X$ admits an
{\em elliptic microspray}
if the $\cO_{\rm an}(X)$-module generated by all
completely integrable holomorphic vector fields  on $X$
is dense in the $\cO_{\rm an}(X)$-module of all
holomorphic vector fields  on $X$ with respect to the
compact-open topology.

We claim that
a Stein manifold $X$ admits an
elliptic microspray if and only if
 $X$ is
holomorphically flexible.
Indeed, admitting an
elliptic microspray implies the holomorphic flexibility,
because
the holomorphic vector fields
on a Stein manifold $X$ span
the tangent space at every point.
As for the converse, we observe that
on a holomorphically
flexible manifold $X$ the sheaf of
germs of holomorphic vector fields
is spanned by the sheaf of
germs of holomorphic vector fields
generated by completely integrable such fields.
By Cartan's Theorem B,  on a Stein  manifold $X$
the corresponding
$\cO_{\rm an}(X)$-modules coincide.
\erems

In the analytic setting, the following analog of Corollary
\ref{locclo} holds.

\blem\label{6.40} If a connected Stein manifold $X$ is
holomorphically flexible over a Stein manifold $B$ then the
relative tangent bundle of $X$ over $B$ is spanned  by a finite
number of completely integrable relative holomorphic vector fields
on $X$. \elem

\bproof In the absolute case i.e., when $B$ is  a point, the
assertion is just that of Lemma 4.1 in \cite{KK}. The proof of
this lemma in \cite{KK} works without changes in the relative case
as well. \eproof

With the same arguments as in the proof of Proposition \ref{6.10}
this implies that a Stein manifold $X$, which is holomorphically
flexible over another Stein manifold $B$,
 admits a fiber dominating spray. Thus we obtain the following result.

\bcor\label{fds} Every connected Stein manifold $X$
holomorphically flexible over another Stein manifold $B$
 admits a fiber dominating spray.
Consequently, the Oka-Grauert-Gromov Principle is valid
for $X\to B$.
\ecor

In particular, the Oka-Grauert principle holds for any
holomorphically flexible connected Stein manifold $X$.

Comparing with the algebraic setting,
in the analytic case
we know little about invariants
of completely integrable holomorphic
vector fields.
This leads to the following question.

\bprob\label{6.60} {\em Does the group $\Aut_{\rm an}(X)$ of
holomorphic automorphisms of a  flexible connected Stein manifold
$X$ act infinitely transitively on $X$?} \eprob

This group is transitive on $X$. Indeed, by the implicit function
theorem every orbit of the group $\Aut_{\rm an}(X)$ is open in
$X$ with respect to the standard Hausdorff topology. On the
other hand, such an orbit is the complement of the union of all
other orbits, thus it is closed. Hence there is only one orbit.

However, the infinite transitivity holds under a stronger
assumption. We need the following notion
from the Andersen-Lempert
theory.

\bdefis\label{6.90} (see \cite{KK},  \cite{Va1})
(i) We say that a complex manifold $X$ has the
{\em density property}
if the Lie algebra generated by all
completely integrable holomorphic vector fields on $X$
is dense in
the Lie algebra of all holomorphic vector fields on
$X$ in the
compact-open topology.

(ii) Similarly, we say that an affine algebraic
manifold $X$ has the
{\em algebraic density
property} if the Lie algebra generated by all completely
integrable algebraic vector fields on
$X$ coincides with  the Lie
algebra of all algebraic vector fields on $X$.
\edefis

An analytic
version of Theorem
\ref{mthm} can be stated as follows
(cf.\ Theorem 5.5 in \cite{For2}).

\bthm\label{varthm} (\cite[2.13]{KK}, \cite{Va1}) If a
connected Stein manifold $X$ of dimension $\ge 2$ has the
density property then the group $\Aut_{\rm an}(X)$ of holomorphic
automorphisms of $X$ acts infinitely transitively\footnote{By
'infinite transitivity` we mean, as before, $m$-transitivity for
all $m\in\N$. Note however that transitivity for arbitrary
discrete subsets does not hold already in $X=\AA^n_\C$, as shows
the famous example of Rosay and Rudin, see e.g., \cite{For2}.} on
$X$. Moreover, for any discrete subset $Z\subseteq X$ and for any
Stein space $Y$ of positive dimension which admits a proper
embedding into $X$, there is another proper embedding $\phi:
Y\hookrightarrow X$ which interpolates $Z$ i.e.,
$Z\subseteq\phi(Y)$. \ethm

We refer the reader to \cite{BF} for a result on interpolation of
a given discrete set of jets of automorphisms by an analytic
automorphism of an affine space, similar to our Theorem
\ref{5.23}.

\subsection{Volume density property}\label{vol}
As usual a holomorphic volume form $\omega$ on a complex manifold
$X$ is a nowhere vanishing top-dimensional holomorphic form on
$X$. We need the following notions.

\bdefis\label{vol0} (i) Given a submersion $X\to B$ of
connected Stein manifolds and a volume form $\omega$ on $X$ we
say that $X$ is {\em holomorphically volume flexible over $B$},
if Definition \ref{6.30}(ii) holds with all relative holomorphic
vector fields considered there being $\omega$-divergence-free. The
latter means that the corresponding phase flow preserves $\omega$.

In the absolute case i.e., $B$ is a point, we simply call the
space $X$ holomorphically volume flexible.

(ii) We say that $X$ has the
{\em volume density property}  if Definition \ref{6.90} holds
with all fields in consideration being
$\omega$-divergence-free.
The {\em algebraic volume density property} is defined likewise.
\edefis

The holomorphic volume flexibility of a Stein manifold $X$ is
equivalent to the existence on $X$ of an elliptic volume
microspray as introduced in \cite{Va}. Lemma \ref{6.40} and
Corollary \ref{fds} admit analogs in this new context. However,
the proofs become now more delicate. We address the interested
reader to \cite{KK, KK1}.

The algebraic volume density property implies the usual volume
density property \cite{KK1}. However, we do not know whether a
holomorphically volume flexible connected Stein manifold has
automatically the volume density property (cf.\ \cite{Va}).

Concerning infinite transitivity, the following theorem is proven
in \cite[2.1-2.2]{KK}.

\bthm\label{ifv} Let $X$ be a connected Stein manifold of
dimension $\ge 2$ equipped with a holomorphic volume form. If $X$
satisfies the holomorphic volume density property, then the
conclusions of Theorem \ref{varthm} hold, with volume preserving
automorphisms. \ethm

Given an algebraic volume form $\omega$ on a smooth affine
algebraic variety $X$, every locally nilpotent vector field on $X$
is automatically $\omega$-divergence-free. Thus the usual
flexibility implies the algebraic volume flexibility. Let us
formulate the following related problem.

\bprob\label{6.100} {\em Let $X$ be a flexible smooth
connected affine algebraic variety over $\C$ equipped with an
algebraic volume form. Does the algebraic volume density property
hold for $X$? } \eprob

We conclude with yet another problem.

\bprob\label{6.110} {\em  Does there exist a flexible exotic
algebraic structure on an affine space that is, a flexible smooth
affine variety over $\C$ diffeomorphic  but not isomorphic to an
affine space $\AA^n_\C$?} \eprob

Notice that for all exotic structures on $\AA^n_\C$ known so far
the Makar-Limanov invariant is non-trivial, whereas for a flexible
such structure, by Proposition \ref{prep} even the field
Makar-Limanov invariant must be trivial.

\brem The preprint of the present paper inspired some further
related results and interesting conjectures, see  \cite{BKK},
\cite{Do2}, and \cite{Pe}. In particular, according to \cite{Pe}
the affine cones over smooth del Pezzo surfaces of degree $\ge 4$
are flexible. In \cite{BKK} a stable birational version of
infinite transitivity is proposed. Conjecture 1.4 in \cite{BKK}
relates this property to unirationality in the sense converse to
that of Proposition \ref{prep}. This should give a
characterization of unirationality versus rational
connectedness. \erem

\providecommand{\bysame}{\leavevmode\hboxto3em{\hrulefill}\thinspace}

\end{document}